\documentclass{siamltex}
\pdfoutput=1
\usepackage{amsmath,amssymb}
\usepackage{array}
\usepackage{xy}
\usepackage{graphicx}
\DeclareGraphicsExtensions{.jpg,.pdf,.png,.jpeg}
\usepackage{epsfig}
\usepackage{mathrsfs}
\usepackage{pstricks}
\usepackage{graphicx}
\usepackage{multirow}
\usepackage{enumerate}
\usepackage{subfig}
\usepackage{color}

\newcommand{\bil}[1]{{#1}}

\newtheorem{remark}{Remark}[section]
\newtheorem{algorithm}{Algorithm}[section]
\newtheorem{example}{Example}[section]

\newcommand{\hF}{\widehat{F}}
\newcommand{\cP}{\mathcal{P}}

\newcommand{\cS}{\mathcal{S}}
\newcommand{\cT}{\mathcal{T}}
\newcommand{\cE}{\mathcal{E}}

\newcommand{\bfR}{\mathbf{R}}
\newcommand{\Ome}{\Omega}
\newcommand{\oOme}{\overline{\Omega}}
\newcommand{\poly}{\mathcal{P}}
\newcommand{\textthm}[1]{\text{\textnormal{#1}}}

\begin{document}

\title{Nonstandard local discontinuous Galerkin methods for  
fully nonlinear second order elliptic and parabolic equations in high dimensions\thanks{This work  
was partially supported by the NSF grant DMS-0710831.}}


\author{Xiaobing Feng\thanks{Department of Mathematics, The University
of Tennessee, Knoxville, TN 37996, U.S.A. ({\tt xfeng@math.utk.edu}) }
\and
Thomas Lewis\thanks{Department of Mathematics and Statistics, The University of 
North Carolina at Greensboro, Greensboro, NC 27410, U.S.A. ({\tt tllewis3@uncg.edu}) 
}
}

\markboth{Xiaobing Feng and Thomas Lewis}{LDG methods for second order fully nonlinear PDEs}

\maketitle
\begin{abstract}
This paper is concerned with developing accurate and efficient 
numerical methods for fully nonlinear second order  
elliptic and parabolic partial differential equations (PDEs) in multiple spatial 
dimensions. It presents a general framework for constructing high order
local discontinuous Galerkin (LDG) methods for approximating
viscosity solutions of these fully nonlinear PDEs. 
The proposed LDG methods are natural extensions of a narrow-stencil finite difference 
framework recently proposed by the authors for approximating viscosity solutions.  
The idea of the methodology is to use multiple approximations of first and second 
order derivatives as a 
way to resolve the potential low regularity of the underlying viscosity solution.  
Consistency and generalized monotonicity properties are proposed 
that ensure the numerical operator approximates the differential operator.
The resulting algebraic system has several linear
equations coupled with only one nonlinear 
equation that is monotone in many of its arguments.  
The structure can be 
explored to design nonlinear solvers. 
This paper also presents and analyzes numerical results for several 
numerical test problems in two dimensions which are used to 
gauge the accuracy and efficiency of the proposed LDG methods.
\end{abstract}




\section{Introduction}\label{sec-1}
In this paper we consider the following general 
fully nonlinear second order elliptic and parabolic PDEs in high dimensions:
\begin{equation}\label{pde_ell}
F[u] := F \left( D^2 u, \nabla u, u, x \right) = 0 , \qquad x \in \Omega 
\end{equation}
and
\begin{equation} \label{pde}
u_t + F \left( D^2 u, \nabla u, u, x, t \right) = 0 , 
\qquad (x,t) \in \Omega_T:=\Omega\times (0,T] 
\end{equation}
which are complemented by appropriate boundary and initial conditions for 
$\Omega\subset\bfR^d (d=2,3)$ a given bounded (possibly convex) domain.
In particular, we are concerned with directly approximating $C^0(\overline{\Omega})$   
(or bounded) solutions of fully nonlinear problems that correspond to the two prototypical 
fully nonlinear operators 
$$F[u] = \text{det } (D^2 u) \qquad\text{and}\qquad  
F[u] = \inf_{\theta \in \Theta} \left( L_\theta u - f_\theta \right), $$ 
where $L_\theta$ is a second order linear elliptic operator with 
$$L_\theta u: = A^\theta : D^2 u + b^\theta \cdot \nabla u + c^\theta u$$  
for $A:B$ the Frobenius 
inner product for matrices $A,B \in \bfR^{d \times d}$.
The first nonlinear operator defines the Monge-Amp\`ere equation, \cite{Pogorelov}, 
and the second nonlinear operator defines the Hamilton-Jacobi-Bellman equation, 
\cite{Fleming,Fleming2}.   It should be noted that some parabolic counterparts of elliptic Monge-Amp\`ere type equations may not have the form of \eqref{pde} (cf.
\cite{Lieberman96}).  Fully nonlinear second order PDEs arise from many 
scientific and engineering fields \cite{FGN12}; 
they are a class of PDEs which are
very difficult to analyze and even more challenging to approximate numerically.

Due to their fully nonlinear structures, fully nonlinear PDEs do not have variational (or weak) formulations 
in general.  
The weak solutions are often defined as {\em viscosity solutions}
(see section \ref{sec-2_1} for the definition).  
The non-variational structure prevents the applicability of standard Galerkin 
type methods such as finite element methods. On the other hand, to approximate very low 
regularity solutions of these PDEs, it is natural to use totally discontinuous piecewise 
polynomial functions (i.e., DG functions) due to 
their flexibility and the larger approximation spaces. As expected, such a method 
must be nonstandard (again) due to the fully nonlinear structure of these PDEs. 
Indeed, a class of nonstandard mixed interior penalty discontinuous Galerkin methods were  
developed by the authors in \cite{Feng_Lewis13} that work well in both $1$-D 
and high dimensions provided that the viscosity solutions belong to 
$C^0 (\overline{\Omega}) \cap H^1(\Omega)$ and 
the polynomial degree is greater than or equal to $1$. Their extensions to 
local discontinuous Galerkin (LDG) methods were done only in the $1$-D case so far.
There were several non-trivial barriers preventing the extensions in the high dimensional case.

The goal of this paper is to generalize the one-dimensional LDG framework and methods of 
\cite{Feng_Lewis12c} to approximate the PDEs \eqref{pde_ell} and \eqref{pde} in 
high dimensions (i.e., $d\geq 2$). Specifically, we shall design and implement a class of  local 
discontinuous Galerkin (LDG) methods which are based on 
a nonstandard mixed formulation of \eqref{pde_ell} and \eqref{pde}.  
Our interest in an LDG approach over the interior-penalty (IP) approach 
found in \cite{Feng_Lewis13} is three-fold.  
The first reason is due to the known increased potential for approximating gradients of 
regular solutions when compared with IPDG methods. 
The second motivation is due to the fact that the LDG approach 
will allow us to form two numerical gradients 
when discretizing fully nonlinear operators that formally involve the gradient of the 
viscosity solution. As already mentioned above, the formulation for the IPDG methods in \cite{Feng_Lewis13} 
assumed the viscosity solutions were in the space 
$C^0 (\overline{\Omega}) \cap H^1(\Omega)$.  
By forming two numerical gradients, the LDG methods can naturally be formulated for 
viscosity solutions in the space $C^0 (\overline{\Omega}) \setminus H^1(\Omega)$.   
Third, as will be seen in the following, the numerical derivatives associated with the LDG approach 
naturally generalize the corresponding difference quotients associated with a 
finite difference (FD) approach. Thus, we can potentially gain further insight into various FD 
methods for fully nonlinear problems by studying their LDG counterparts while also having a
stronger theoretical foundation for the LDG methods proposed in this paper.  

The main difficulty addressed in this paper is how to extend the one-dimensional framework 
of  \cite{Feng_Lewis12c} to the high-dimensional setting.  
First, we will need to design a consistent way for forming multiple discrete 
gradient and Hessian approximations.  
To this end, we will utilize the conventions introduced in \cite{Feng_Lewis_Neilan13} 
where a finite element DG numerical calculus was developed based upon a discontinuous Galerkin 
methodology and choosing various fluxes to characterize various numerical derivative operators.  
To extend ideas to the high-dimensional setting we will discretize partial derivatives directly 
as a way to define various gradient approximations.    
We will need to introduce nonstandard trace operators 
that are consistent with the idea that each partial derivative is treated independently.  
Second, we will extend the framework to second order problems where the fully nonlinear 
differential operator also involves the gradient operator, as represented by the general 
problems \eqref{pde_ell} and \eqref{pde}.  
Third,  on noting that the LDG formulation will introduce a large set of auxiliary equations, 
we will explore various solver techniques and the potential for variable reduction to reduce 
the computational cost. 

We note that typically a DG formulation for a fully nonlinear problem is based upon a 
semi-Lagrangian approach or strong structure assumptions that guarantee a monotonicity 
property of the scheme (see \cite{Debrabant_Jakobsen13,Feng_Jensen,Jensen_Smears,NSZ16,NNZ17,Salgado_Zhang16,Smears_Suli} 
and the review article \cite{FGN12}). 
As such, the methods are limited to piecewise linear basis functions.  
Inspired by the work of Yan and Osher in \cite{Yan_Osher11}, 
we seek to formulate DG methods that allow the use of high order polynomials   
and can achieve high-order accuracy. The methods proposed in this paper extend 
the narrow-stencil FD approach in \cite{FKL11,Lewis_Dissertation,FL-FDnD} to high-order
and to unstructured triangular meshes. 
As with the LDG methods for Hamilton-Jacobi equations in \cite{Yan_Osher11}, 
the only analytic convergence result for the proposed LDG methods corresponds to 
choosing piecewise constant basis functions.  
In this special case, the proposed LDG methods reduce to the FD methods of \cite{FL-FDnD} 
and may lead to developing new high-order narrow-stencil FD methods.  
Moreover, we are able to extablish a link between the proposed LDG 
methods and the vanishing moment method of Feng and Neilan \cite{Feng_Neilan11}. 
Heuristically such a link also helps to justify the proposed 
LDG methods for fully nonlinear second order PDEs in the same way
the link to the vanishing viscosity method motivates the LDG methods of \cite{Yan_Osher11} 
for Hamilton-Jacobi equations. 

The remainder of this paper is organized as follows. 
In section \ref{sec-2_1} we introduce some background for the viscosity solution notion.  
In section \ref{sec-2_2} we define key concepts of consistency and generalized monotonicity 
for numerical operators that will serve as the foundation of the proposed LDG framework.  
We also introduce the numerical operators that will be used in the design of our methods.  
The proposed LDG formulation for the nonlinear elliptic equation \eqref{pde_ell} 
is presented in section \ref{sec-3}. We use two main ideas in the formulation:
the numerical viscosity borrowed from the discretization of first-order Hamilton-Jacobi equations 
and a novel concept of numerical moments.  
We also discuss various techniques for solving the resulting nonlinear (large) algebraic systems.  
In section \ref{sec-4} we consider both explicit and 
implicit in time fully discrete LDG methods for the fully nonlinear 
parabolic equation \eqref{pde} based on the method of lines approach. 
In section \ref{sec-5} we present many numerical experiments for the proposed LDG methods.
These numerical experiments verify the accuracy and demonstrate the efficiency of the new methods.
The experiments also explore the role of the numerical moment in the formulation.
Lastly, in section \ref{sec-6}, we provide some concluding remarks.  

\section{Preliminaries} \label{sec-2_1}

We first recall the viscosity solution concept for fully nonlinear second order problems.  
For a bounded open domain $\Ome\subset\mathbf{R}^d$, let $B(\Ome)$, 
$USC(\Ome)$, and $LSC(\Ome)$ denote, respectively, the spaces of bounded,
upper semi-continuous, and lower semi-continuous functions on $\Ome$.
For any $v\in B(\Ome)$, we define
\[
v^*(x):=\limsup_{y\to x} v(y) \qquad\mbox{and}\qquad
v_*(x):=\liminf_{y\to x} v(y). 
\]
Then, $v^*\in USC(\Ome)$ and $v_*\in LSC(\Ome)$, and they are called
{\em the upper and lower semicontinuous envelopes} of $v$, respectively.

Given a function $F: \cS^{d\times d}\times\mathbf{R}^d\times 
\mathbf{R}\times \oOme \to \mathbf{R}$, where $\cS^{d\times d}$ denotes the set
of $d\times d$ symmetric real matrices, the general second order
fully nonlinear PDE takes the form
\begin{align}\label{e2.1}
F(D^2u,\nabla u, u, x) = 0 \qquad\mbox{in } \oOme.
\end{align}
Note that here we have used the convention of writing the boundary condition as a
discontinuity of the PDE (cf. \cite[p.274]{Barles_Souganidis91}).

The following two definitions can be found in \cite{Gilbarg_Trudinger01,Caffarelli_Cabre95,Barles_Souganidis91}.

\begin{definition}\label{def2.1}
Equation \eqref{e2.1} is said to be {\em elliptic} if for all
$(\mathbf{q},\lambda,x)\in \mathbf{R}^d\times \mathbf{R}\times \oOme$ there holds
\begin{align}\label{e2.2}
F(A, \mathbf{q}, \lambda, x) \leq F(B, \mathbf{q}, \lambda, x) \qquad\forall 
A,B\in \cS^{d\times d},\, A\geq B, 
\end{align}
where $A\geq B$ means that $A-B$ is a nonnegative definite matrix.
\end{definition}
Equation \eqref{e2.1} is said to be {\em proper elliptic} if for all
$(\mathbf{q},x)\in \mathbf{R}^d \times \oOme$ there holds
\begin{align}\label{e2.2b}
F(A, \mathbf{q}, a, x) \leq F(B, \mathbf{q}, b, x) \qquad\forall 
A,B\in \cS^{d\times d},\, A\geq B, \; a,b \in \mathbf{R},\,a\leq b . 
\end{align}
We note that when $F$ is differentiable, ellipticity
can also be defined by requiring that the matrix $\frac{\partial F}{\partial A}$
is negative semi-definite (cf. \cite[p. 441]{Gilbarg_Trudinger01}).

\begin{definition}\label{def2.2}
A function $u\in B(\Ome)$ is called a viscosity subsolution (resp.
supersolution) of \eqref{e2.1} if, for all $\varphi\in C^2(\oOme)$,
if $u^*-\varphi$ (resp. $u_*-\varphi$) has a local maximum
(resp. minimum) at $x_0\in \oOme$, then we have
\[
F_*(D^2\varphi(x_0),\nabla \varphi(x_0), u^*(x_0), x_0) \leq 0 
\]
(resp. $F^*(D^2\varphi(x_0),\nabla \varphi(x_0), u_*(x_0), x_0) \geq 0$).
The function $u$ is said to be a viscosity solution of \eqref{e2.1}
if it is simultaneously a viscosity subsolution and a viscosity
supersolution of \eqref{e2.1}.
\end{definition}

\begin{remark}\label{rem2.1}
It can be proved that it is sufficient only to consider $\varphi \in\poly_2$,
the space of all {\em quadratic polynomials}, in Definition \ref{def2.2} (see
\cite[page 20]{Caffarelli_Cabre95}).
\end{remark}

\begin{definition}\label{comparison}
Problem \eqref{e2.1} is said to satisfy a {\em comparison 
principle} if the following statement holds. For any upper semi-continuous 
function $u$ and lower semi-continuous function $v$ on $\overline{\Omega}$,  
if $u$ is a viscosity subsolution and $v$ is a viscosity supersolution 
of \eqref{e2.1}, then $u\leq v$ on $\overline{\Omega}$.
\end{definition}

We remark that if $F$ and $u$ are continuous, then the upper and lower $*$
indices can be removed in Definition \ref{def2.2}. The definition
of ellipticity implies that the differential operator $F$
must be non-increasing in its first argument in order to be
elliptic. It turns out that ellipticity and a comparison principle provide sufficient
conditions for equation \eqref{e2.1} to fulfill a maximum principle
(cf. \cite{Gilbarg_Trudinger01,Caffarelli_Cabre95}).
It is clear from the above definition that viscosity solutions
in general do not satisfy the underlying PDEs in a tangible sense, and
the concept of viscosity solutions is {\em nonvariational}. Such
a solution is not defined through integration by parts against arbitrary test
functions; hence, it does not satisfy an integral identity. 
The non-variational nature of viscosity
solutions is the main obstacle that prevents the direct construction
of Galerkin-type methods.


\section{A generalized monotone nonstandard LDG framework} \label{sec-2_2}

Our methodology for directly approximating viscosity solutions of second-order 
fully nonlinear PDEs is based on several motivational ideas which we explain below.  
Since integration by parts cannot be performed on equation \eqref{pde_ell}, 
{\em the first key idea} is to introduce the auxiliary variables $P:=D^2 u$ and $q:= \nabla u$ 
and rewrite the original fully nonlinear PDE as a system of PDEs:  
\begin{subequations} \label{mixed_1}
\begin{eqnarray}
F(p,q,u,x) &=0, \label{mixed_1F} \\
q- \nabla u &=0, \\
P-\nabla q &=0. 
\end{eqnarray}
\end{subequations}
To address the fact that $\nabla u$ and $D^2 u$ may not exist for a viscosity 
solution $u\in C^0(\overline{\Omega})$, 
{\em the second key idea} 
is to formally replace $q:=\nabla u$ by two possible values of 
$\nabla u$, namely, the left and right (possibly infinite) limits, and $P:= \nabla q$ by two possible
values for each possible $q$, namely, the left and right (possibly infinite) limits.  
Thus, we have the auxiliary variables
$q^-, q^+ : \Omega \to \bfR^d$ 
and $P^{- -}, P^{- +}, P^{+ -}, P^{+ +} : \Omega \to \bfR^{d \times d}$
such that
\begin{subequations}\label{mixed_2}
\begin{align} 
\left[ q^-(x) \right]_i &=  \lim_{\sigma \to 0^+} \left[ \nabla u(x - \sigma \mathbf{e}_i) \right]_i ,
 \label{qdn} \\
\left[ q^+(x) \right]_i &= \lim_{\sigma \to 0^+} \left[ \nabla u(x + \sigma \mathbf{e}_i) \right]_i, \label{qup} \\
\left[ P^{- -}(x) \right]_{ij} &= \lim_{\sigma \to 0^+} \left[ \nabla q^-(x - \sigma \mathbf{e}_j) \right]_{ij}, \label{pdndn} \\
\left[ P^{- +}(x) \right]_{ij} &= \lim_{\sigma \to 0^+} \left[ \nabla q^-(x + \sigma \mathbf{e}_j) \right]_{ij}, \label{pdnup} \\
\left[ P^{+ -}(x) \right]_{ij} &=\lim_{\sigma \to 0^+} \left[ \nabla q^+(x - \sigma \mathbf{e}_j) \right]_{ij}, \label{pupdn} \\
\left[ P^{+ +}(x) \right]_{ij} &=\lim_{\sigma \to 0^+} \left[ \nabla q^+(x + \sigma \mathbf{e}_j) \right]_{ij} \label{pupup} 
\end{align}
\end{subequations} 
for all $i,j \in \{1,2,\ldots,d\}$, 
where $\mathbf{e}_i$ denotes the $i$th canonical basis vector for $\mathbf{R}^d$.  
{\em The third key idea} is to replace \eqref{mixed_1F} by 
\begin{eqnarray} \label{F_hat}
\hF(P^{+ +}, P^{+ -}, P^{- +}, P^{- -}, q^+, q^-, u, x) = 0,
\end{eqnarray} 
where $\hF$, which is called a {\em numerical operator}, should be some
well-chosen approximation to $F$ that incorporates the multiple gradient and Hessian variables. 

The next step is to address the key issue about  what criterion or properties  ``good" numerical 
operators $\hF$ should satisfy.  A large part of our framework revolves around describing sufficient conditions on the choice of numerical operators, as reflected in the following definitions that 
generalize the one-dimensional definitions given in \cite{Feng_Lewis12c}.  

\smallskip
\begin{definition}\label{def2.3}
\begin{itemize}
\item[{\rm (i)}] A function 
$\hF: \bigl( \bfR^{d \times d} \bigr)^4 \times \bigl( \bfR^d \bigr)^2 \times \bfR \times \Omega \to \bfR$ 
is called  a {\em numerical operator}. 
\item[{\rm (ii)}] Let $P \in \overline{\bfR}^{d \times d}$, $q \in \overline{\bfR}^d$, $v \in \bfR$, and 
$x \in \overline{\Omega}$.  
A numerical operator $\hF$ is said to be {\em consistent} (with 
the differential operator $F$) if $\hF$ satisfies
\begin{align}\label{A1a}
\liminf_{P^{\mu \nu} \to P; \mu, \nu = -, + \atop q^\pm \to q, \lambda \to v, \xi \to x} 
\hF(P^{+ +}, P^{+ -}, P^{- +}, P^{- -}, q^+, q^-,\lambda, \xi) \geq F_*(P,q,v,x),\\
\limsup_{P^{\mu \nu} \to P; \mu, \nu = -, + \atop q^\pm \to q, \lambda \to v, \xi \to x} 
\hF(P^{+ +}, P^{+ -}, P^{- +}, P^{- -}, q^+, q^-,\lambda, \xi) \leq F_*(P,q,v,x), \label{A1b}
\end{align}
where $F_*$ and $F^*$ denote, respectively, the lower and the upper
semi-continuous envelopes of $F$. Thus, we have 
\begin{align*}
	F_*(P,q,v,x) & := 
		\liminf_{\widetilde{P} \to P, \widetilde{q} \to q, \atop \widetilde{v} \to v, \widetilde{x} \to x} 
			F \bigl( \widetilde{P}, \widetilde{q}, \widetilde{v}, \widetilde{x} \bigr), \\ 
	F^*(P,q,v,x) & := 
		\limsup_{\widetilde{P} \to P, \widetilde{q} \to q, \atop \widetilde{v} \to v, \widetilde{x} \to x} 
			F \bigl( \widetilde{P}, \widetilde{q}, \widetilde{v}, \widetilde{x} \bigr),
\end{align*}
where $\widetilde{P} \in \bfR^{d \times d}$, $\widetilde{q} \in \bfR^d$, $\widetilde{v} \in \bfR$, 
and $\widetilde{x} \in \Omega$.  
Note, when $F$ and $\widehat{F}$ are continuous, the above definition can be simplified to
\begin{equation}\label{A1c}
\hF(P, P, P, P ,q, q,v, x) = F(P,q,v,x) . 
\end{equation}

\item[{\rm (iii)}] A numerical operator $\hF$ is said to be {\em g-monotone} 
if for all $x \in \Omega$, there holds 
$\hF(P^{+ +} , P^{+ -} , P^{- +} , P^{- -} , q^+, q^-, v, x)$ 
is monotone increasing in $P^{+ +}$, $P^{- -}$, $q^-$, and $v$ 
and monotone decreasing in $P^{+ -}$, $P^{- +}$, and $q^+$.
More precisely, the numerical operator $\hF$ is g-monotone 
if for all $P^{\mu \, \nu} \in \mathbf{R}^{d \times d}$ and 
$q^\mu \in \mathbf{R}^d$, $\mu, \nu \in \{+,-\}$, 
for all $v \in \mathbf{R}$, and for all $x \in \Omega$, there holds
\begin{align*}
\hF \bigl( A , P^{+ -} , P^{- +} , P^{- -} , q^+ , q^- , v , x \bigr) 
	& \leq \hF \bigl( B , P^{+ -} , P^{- +} , P^{- -} , q^+ , q^- , v , x \bigr) , \\
\hF \bigl( P^{+ +} , A , P^{- +} , P^{- -} , q^+ , q^- , v , x \bigr) 
	& \geq \hF \bigl( P^{+ +} , B , P^{- +} , P^{- -} , q^+ , q^- , v , x \bigr) , \\
\hF \bigl( P^{+ +} , P^{+ -} , A , P^{- -} , q^+ , q^- , v , x \bigr) 
	& \geq \hF \bigl( P^{+ +} , P^{+ -} , B , P^{- -} , q^+ , q^- , v , x \bigr) , \\
\hF \bigl( P^{+ +} , P^{+ -} , P^{- +} , A , q^+ , q^- , v , x \bigr) 
	& \leq \hF \bigl( P^{+ +} , P^{+ -} , P^{- +} , B , q^+ , q^- , v , x \bigr) , 
\end{align*}
for all $A,B \in \cS^{d\times d}$ such that $A\preceq B$, where $A \preceq B$ means that
$B-A$ has all nonnegative components, 
\begin{align*}
\hF \bigl( P^{+ +} , P^{+ -} , P^{- +} , P^{- -} , a , q^- , v , x \bigr) 
	& \geq \hF \bigl( P^{+ +} , P^{+ -} , P^{- +} , P^{- -} , b , q^- , v , x \bigr) , \\
\hF \bigl( P^{+ +} , P^{+ -} , P^{- +} , P^{- -} , q^+ , a , v , x \bigr) 
	& \leq \hF \bigl( P^{+ +} , P^{+ -} , P^{- +} , P^{- -} , q^+ , b , v , x \bigr) , 
\end{align*}
for all $a,b \in \mathbf{R}^d$ such that 
$a_i \leq b_i$ for all $i = 1, 2, \ldots, d$, 
and 
\[
\hF \bigl( P^{+ +} , P^{+ -} , P^{- +} , P^{- -} , q^+ , q^- , a , x \bigr) 
	\leq \hF \bigl( P^{+ +} , P^{+ -} , P^{- +} , P^{- -} , q^+ , q^- , b , x \bigr) 
\] 
for all $a,b \in \mathbf{R}$ such that $a \leq b$.  

The condition can be summarized by 
$\hF(\uparrow,\downarrow,\downarrow,\uparrow, \downarrow, \uparrow, \uparrow, x)$, 
where the monotonicity with respect to the matrix entries is enforced component-wise.  
\end{itemize}

\end{definition}

The final concern for the framework is how to design numerical operators that 
are both consistent and g-monotone.  
Inspired by Lax-Friedrichs numerical Hamiltonians used for Hamilton-Jacobi equations \cite{Shu07}, 
we propose the following Lax-Friedrichs-like numerical operator:
\begin{align}\label{LF1}
&\hF(P^{+ +}, P^{+ -}, P^{- +}, P^{- -}, q^+, q^-,\lambda, \xi)
:= F \bigg(\frac{P^{- +} + P^{+ -}}{2},\frac{q^- + q^+}{2},\lambda,\xi \bigg) \\
\nonumber & \hskip 1.5in 
- \beta \cdot \bigl( q^- - q^+ \bigr) 
+ \alpha : \bigl(P^{+ +} - P^{+ -} - P^{- +} + P^{- -} \bigr),
\end{align}
where $\alpha \in \bfR^{d \times d}$ is an undetermined positive semi-definite matrix and 
$\beta \in \bfR^d$ is an undetermined nonnegative vector.  $A:B$ stands for the Frobenius 
inner product 
for matrices $A,B \in \bfR^{d \times d}$.
The second to last term $\beta \cdot ( q^- - q^+ )$ is referred to as the numerical viscosity 
and is directly borrowed from Lax-Friedrichs numerical Hamiltonians, 
and the last term $\alpha : ( P^{+ +} - P^{+ -} - P^{- +} + P^{- -} )$ is referred to as the {\em numerical moment}.  
It is trivial to verify that $\hF$ is consistent with $F$ when $F$ is continuous.
By choosing $\alpha$ and $\beta$ correctly, we can also ensure g-monotonicity.
In practice, we typically choose $\beta = b \vec{1}$ and $\alpha = a_1 I + a_2 \mathbf{1}$ 
for sufficiently large positive constants $a_1$, $a_2$, and $b$, 
where $\mathbf{1}$ is the vector/matrix with all entries equal to one and $I$ is the identity matrix.  
We note that the g-monotonicity condition can be realized for $a_2$ sufficiently large and $a_1 = 0$.  
By also choosing $a_1$ large, we can additionally enforce the g-monotonicity condition using the partial 
order based on SPD matrices.  

\begin{remark} \ 
\begin{enumerate}[(a)]
\item 
Due to the definition of ellipticity for $F$, the g-monotonicity constraints on $\hF$ with respect
to $P^{- +}_{ii}$ and $P^{+ -}_{ii}$ are natural.  
Consistency is used to pass to a single matrix argument and ellipticity is used 
to guarantee the correct monotonicity with respect to the partial ordering induced by 
SPD matrices.  
\item 
By choosing the numerical viscosity and the numerical moment correctly, 
the numerical operator $\hF$ will behave
like a strongly elliptic operator even if the PDE operator $F$ is a degenerate elliptic operator.
The consistency assumption then guarantees that the numerical operator is still a reasonable
approximation for the PDE operator.
\item 
When  $F$ is differentiable, while it may not be possible to globally bound 
$\frac{\partial F}{\partial \nabla u}$ and 
$\frac{\partial F}{\partial D^2 u}$, it may be sufficient
to choose values for $\beta$ and $\alpha$ such that the g-monotonicity property is preserved locally
over each iteration of the nonlinear solver for a given initial guess.  
The same remark holds if $F$ is locally Lipschitz.  
\end{enumerate}
\end{remark}

\section{Formulation of nonstandard LDG methods for elliptic PDEs} \label{sec-3}

We now formulate our nonstandard LDG methods for approximating viscosity solutions 
of fully nonlinear elliptic PDEs which are based on 
the mixed formulation \eqref{mixed_2} and \eqref{F_hat}.  
We also provide a detailed explanation of how to treat the boundary traces 
in the formulation.  
Lastly we use the DG formulation to better understand the numerical viscosity and numerical 
moment appearing in our Lax-Friedrichs-like numerical operator 
and explore two algorithms for solving the resulting nonlinear algebraic systems.  

\subsection{DG Notation} \label{sec-2_3}

To formulate our LDG methods, we need to introduce some notation and conventions 
 which are standard and can be found in \cite{Feng_Lewis_Neilan13}.   
Let $\Omega$ be a polygonal domain and $\mathcal{T}_h$ 
denote a locally quasi-uniform and shape-regular partition 
of $\Omega$
with $h = \max_{K \in \cT_h} (\text{diam} K)$.
We introduce the broken $H^1$-space and broken $C^0$-space
\[
	H^1(\cT_h) := \prod_{K \in \cT_h} H^1(K) , \qquad 
	C^0 (\cT_h) := \prod_{K \in \cT_h} C^0 ( \overline{K} )  
\]
and the broken $L^2$-inner product
\[
(v ,w)_{\mathcal{T}_h} := \sum_{K \in \cT_h} \int_{K} v w\, dx \qquad 
	\forall v,w \in L^2(\cT_h) .
\]
Let $\cE_h^I$ denote the set of all interior faces/edges of $\cT_h$, 
$\cE_h^B$ denote the set of all boundary faces/edges of $\cT_h$, 
and $\cE_h := \cE_h^I \cup \cE_h^B$.
Then, for a set $\cS_h \subset \cE_h$, we define the broken $L^2$-inner product over $\cS_h$ by 
\[
\langle v ,w \rangle_{\cS_h} := \sum_{e \in \cS_h} \int_{e} v \, w\, ds \qquad 
	\forall v,w \in L^2(\cS_h) .
\]
For a fixed integer $r \geq 0$, we define the standard DG finite element space
$V^h \subset H^1 (\cT_h) \subset L^2(\Omega)$ by
\[
V^h := \prod_{K \in \mathcal{T}_h} \poly_{r} (K),
\]
where $ \poly_{r} (K)$ denotes the set of all polynomials on $K$ with
degree not exceeding $r$.  

For $K, K'\in \cT_h$, let $e=\partial K\cap \partial K' \in \cE^I_h$. Without a loss of
generality, we assume that the global labeling number of $K$ is smaller than 
that of $K'$ and define the following (standard) jump and average notations: 
\begin{equation} \label{DG_jump_avg}
[v] := v|_K-v|_{K'} , \qquad 
\{ v \} := \frac{v|_K + v|_{K'}}{2} 
\end{equation}
for any $v\in H^m(\cT_h)$. We also define $n_e:=n_K=-n_{K'}$ as the normal vector to $e$. 
When $e \in \cE^B_h$, $n_e$ denotes the unit outward normal for the underlying boundary simplex.
We note that the function values defined on $\cE^B_h$ will be handled in a nonstandard way in
 our LDG methods by allowing the boundary function values to depend on 
the degree of the polynomial basis $r$. However, when $r \geq 1$,  the boundary function 
values can be treated in a more standard way as in \cite{Feng_Lewis_Neilan13}.

\subsection{Formulation of LDG methods} \label{sec-3_2} 

We now present an element-wise formulation for our LDG methods.
First we introduce some local definitions. For any $e\in \cE_h^I$ with 
 $e = \partial K \cap \partial K^\prime$  for some $K,K'\in \cT_h$ and for any 
 $v \in V^h$, let $v( x^I)$ denote the value of $v(x)$ on $\partial K$
 from the interior of the element $K$ and $v(x^E)$ denote the
 value of $v(x)$ on $\partial K$ from the interior of the element $K^\prime$.
Using these limit definitions, we then define the local boundary flux operators:
$T^+, T^- : \mathcal{P}_r (K) \to \left( \prod_{e \subset \partial K}\mathcal{P}_r (e) \right)^d$ by
\begin{subequations}\label{DG_discrete_traces}
\begin{align}
	T_i^- (v_h)(x) & := 
		\begin{cases}
			v_h ( x^I ) , \quad & \text{if } n_i(x) > 0 , \\
			v_h ( x^E ) , \quad & \text{if } n_i(x) < 0 , \\ 
			\{ v_h(x) \} , \quad & \text{if } n_i(x) = 0 , 
		\end{cases} \\
	T_i^+ (v_h)(x) & := 
		\begin{cases}
			v_h ( x^E) , \quad & \text{if } n_i(x) > 0 , \\
			v_h ( x^I ) , \quad & \text{if } n_i(x) < 0 , \\ 
			\{ v_h(x) \} , \quad & \text{if } n_i(x) = 0 
		\end{cases}
\end{align}
\end{subequations}
for all $i \in \{ 1, 2 , \ldots , d\}$, $x \in e$, and $v_h \in V^h$.
The definition of $T_i^\pm (v)$ for $v \in V^h$ on each $e \in \cE^B_h$ 
will be delayed to section~\ref{sec-3_1}.  
Observe that, for $e \in \cE^I_h$, we can also rewrite the labelling-dependent trace operators as
\begin{equation} \label{interior_trace}
	T_i^\pm(v_h) = \big\{ v_h \big\} \mp \frac12 \text{sgn} (n_e^{(i)}) \big[ v_h \big]  \quad
\mbox{where} \quad
	\text{sgn}(y) = \begin{cases}
	1 & \text{if } y > 0 , \\ 
	-1 & \text{if } y < 0 , \\ 
	0 & \text{if } y = 0 
	\end{cases} 
\end{equation}
for all $y \in \bfR$, where $n_e^{(i)}$ denotes the $i$-th component of $n_e$ (
the unit outward normal to $e$).  
Note that the trace operators are nonstandard in that their values depend on the individual 
components of the edge normal $n_e$. The standard definition assigns a single-value (called
a numerical flux) based on the edge normal vector as a whole.  

We are now ready to formulate our LDG methods for system \eqref{mixed_2}--\eqref{F_hat}.
First, we approximate the (fully) nonlinear equation \eqref{F_hat} 
by its broken $L^2$-projection into $V^h$, namely,
\begin{equation}\label{pde_ell_weak}
\bil{a}_0 \bigl(u_h , q^+_h, q^-_h, P^{+ +}_h, P^{+ -}_h, P^{- +}_h, P^{- -}_h; \phi_{0h} \bigr)  = 0 
\qquad \forall \phi_{0h} \in V^h, 
\end{equation}
where
\begin{align*}
& \bil{a}_0 (u , q^+, q^-, P^{+ +}, P^{+ -}, P^{- +}, P^{- -}; \phi_0) \\
& \hskip 1in 
= \bigl(\hF(P^{+ +}, P^{+ -}, P^{- +}, P^{- -} ,q^+, q^-,u,\cdot),\phi_0 \bigr)_{\mathcal{T}_h}. 
\end{align*}

Next, we discretize the six {\em linear} equations in \eqref{mixed_2} locally
with respect to each component using the integration by parts formula: 
\begin{equation} \label{DG_greenes}
	\int_S v_{x_i} \, \varphi \, dx
	= \int_{\partial S} v \, \varphi \, n_i \, ds 
		- \int_S v \, \varphi_{x_i} \, dx \qquad \forall \varphi \in C^1(S)
\end{equation}
for $i = 1, 2, \ldots, d$.
Thus, the above formula yields an integral characterization for the partial derivative 
$v_{x_i}$ on the set $S$ for all $v \in H^1 (S)$.
Using the preceding identity, we define our gradient approximations 
$q_h^\mu \in ( V^h)^d$, $\mu \in \{+ , - \}$, by
\begin{equation} \label{DG_local_q}
	\int_K q^\mu_i \, \phi_i^\mu \, dx + \int_K u \, ( \phi_i^\mu)_{x_i} \, dx
		= \int_{\partial K} T_i^\mu (u) \, n_i \, \phi_i^\mu ( x^I ) \, ds \quad
		\forall \phi_i^\mu \in V^h
\end{equation}
for $i = 1, 2, \ldots, d$, $\mu = +, -$.  

Similarly, we define our Hessian approximations 
$P_h^{\mu \, \nu} \in ( V^h )^{d \times d}$, $\mu , \nu \in \{ + , - \}$, by 
\begin{equation} \label{DG_local_p}
	\int_K P^{\mu \, \nu}_{i, j} \, \psi^{\mu \, \nu}_{i, j} \, dx 
		+ \int_K q^\mu_{i} \, ( \psi^{\mu \, \nu}_{i, j} )_{x_j} \, dx
	= \int_{\partial K} T_{j}^{\nu} (q^\mu_i) \, n_j \, \psi^{\mu \, \nu}_{i, j} (x^I) \, ds
\end{equation}
for all $\psi^{\mu\,\nu}_{i,j}\in V^h$ and  $i, j = 1, 2, \dots, d$, $\mu, \nu = + , -$.

Thus, in order to approximate the viscosity solution $u$ for the fully nonlinear PDE 
\eqref{pde_ell} paired with a Dirichlet boundary condition 
\begin{equation}
u= g \qquad\mbox{ on } \partial \Omega \label{bc_ell} 
\end{equation}
for a given function $g \in C^0 (\partial \Omega)$, we seek functions
$u_h$ $\in V^h$; 
$q^+_h$,$ q^-_h$ $\in (V^h)^d$; and 
$P^{+ +}_h$, $P^{+ -}_h,$ $P^{- +}_h$, $P^{- -}_h$ $\in (V^h)^{d \times d}$
such that equation \eqref{pde_ell_weak} holds as well as equations 
\eqref{DG_local_q} and \eqref{DG_local_p} for all $K \in \cT_h$, 
where $u_h$ forms the approximation for $u$. We note that the implementation of 
the Dirichlet boundary condition into the definition of the boundary flux/trace operator
in \eqref{DG_local_q} and \eqref{DG_local_p} will be described in section~\ref{sec-3_1}.    

By summing the definitions of $q_h^\pm$ and $P^{\mu, \nu}_h$ over $\mathcal{T}_h$ and using 
\eqref{interior_trace}, we obtain  the following global (labeling-dependent) formulations for 
the proposed LDG methods:
\begin{subequations}\label{global_aux} 
\begin{align} 
	\big( q_i^\mu , \varphi_i^\mu \big)_{\mathcal{T}_h} + a_i^\mu \big( u_h , \varphi_i^\mu \big) & = 0 
		\qquad \forall \varphi_i^\mu \in V^h , \label{q_bil} \\ 
	\big( P^{\mu \nu}_{i j} , \psi^{\mu \nu}_{i j} \big)_{\mathcal{T}_h} 
		+ a_j^\nu \big( q_i^\mu , \psi^{\mu \nu}_{i j} \big) & = 0 
		\qquad \forall \psi^{\mu \nu}_{i j} \in V^h \label{p_bil}
\end{align}
\end{subequations} 
for $i,j = 1, 2, \ldots, d$ and $\mu, \nu = -, +$, where 
\begin{equation}\label{a_bil}
	a_i^\pm \big( v , \phi \big) := 
		\big( v , \phi_{x_i} \big)_{\mathcal{T}_h} 
		- \Bigl\langle \{v\} \mp \frac12 \text{sgn} (n_e^{(i)}) [v] , [\phi] n_e^{(i)} \Bigr\rangle_{\cE^I_h} 
		- \bigr\langle T_i^\pm(v) , \phi \, n_i \big\rangle_{\cE^B_h} 
\end{equation}
for all $v, \phi \in V^h$.  
Then, the proposed LDG methods correspond to solving the global formulation 
\eqref{pde_ell_weak} and \eqref{global_aux}.  

\begin{remark} 
Since the approximations are piecewise totally discontinuous polynomials, 
the sided limits in \eqref{mixed_2} only need to be enforced along the faces/edges. 
By \cite{Feng_Lewis_Neilan13}, we know  that the proposed auxiliary variables 
provide proper meanings for the limits in \eqref{mixed_2} since the various derivative 
approximations coincide with the $L^2$ projections of distributional derivatives 
onto $V^h$ with variable strengths on the interior faces/edges depending on the 
choices of the traces, where the traces are chosen such that the sided
 limits in \eqref{mixed_2} are consistent.  
\end{remark} 

\subsection{Numerical boundary fluxes} \label{sec-3_1}

In this section, we extend the definition for the boundary flux operators, 
given by \eqref{DG_discrete_traces},  to the set $\cE^B_h$.  
To this end, we will introduce a set of constraint equations that express all exterior limits 
in terms of interior limits and known data. 
The Dirichlet boundary data will serve as an exterior constraint on the sought-after numerical solution.  
We will consider two cases based on whether the order of the DG space $V^h$ is zero or nonzero, i.e.,
 $r=0$ or $r \geq 1$. When $r \geq 1$, we will enforce a ``continuity" assumption across the boundary 
 $\partial \Omega$, and when $r=0$, we will prescribe an alternative approach that will more closely
 resemble the introduction of ``ghost values" commonly used in FD methods.  

Prior to introducing the constraint equations, 
we specify a convention to be used for all boundary faces/edges.   
Let $K \in \cT_h$ be a boundary simplex, and let $e \in \cE^B_h$ such that $e \subset \partial K$. 
Suppose $v_h \in V^h$ such that $v_h$ is supported on $K$.  
Then, we define $v_h(x) := v_h(x^I)$ for all $x \in e$. 

We first consider $r \geq 1$, in which case we make the ``continuity" assumption 
\begin{equation}
	v_h(x^E) = v_h(x) 
\end{equation}
for all $x \in e$ and $v_h \in V^h$ such that $e \in \cE^B_h$.
Since problem \eqref{pde_ell} and \eqref{bc_ell} does not provide a Neumman 
boundary data, we simply treat $q_i^\pm(x)$ as an unknown 
for all $i=1,2,\ldots,d$ and $x \in e$ with $e \in \cE^B_h$.
Alternatively, when defining the boundary flux values for $u_h$, we use 
the Dirichlet boundary condition given by \eqref{bc_ell}.  
Thus, for $r \geq 1$, we wish to impose
\[ 
	u_h \left( x \right) = g(x)  
\]
for all $x \in \partial \Omega$.
However, $g$ may not be a polynomial of degree $r$.  
Thus, we enforce this condition weakly by imposing the following constraint equations:
\begin{equation} \label{DG_u_bc}
	\sum_{i=1}^d \bigl\langle u_h(x) , \varphi_h(x) n_i \bigr\rangle_{\cE^B_h} 
		= \sum_{i=1}^d \bigl\langle g(x) , \varphi_h(x) n_i \bigr\rangle_{\cE^B_h} 
		\qquad \forall \varphi_h \in V^h,
\end{equation} 
where $n$ denotes the unit outward normal vector along $\partial \Omega$.  
Observe that when a boundary simplex has more than one face/edge in $\cE^B_h$, 
we are treating all of the boundary simplex's faces/edges in $\cE^B_h$ 
as a single ($d$-1)-dimensional surface.  



We now consider the case $r = 0$.
Extending the definition for the boundary flux operators, given by \eqref{DG_discrete_traces},  
to the set $\cE^B_h$ is less straightforward in this case.
We can see this by observing the fact that 
when fixing the interior limit of a boundary value on a boundary simplex, 
we actually fix the function value on the entire simplex.
Thus, strictly enforcing a Dirichlet boundary condition for $u_h$ may result in a boundary 
layer with respect to the overall approximation error when measured in low-order norms such as the 
$L^\infty$- or $L^2$-norm.  
Our goal is to prescribe boundary flux values in a way that results in a potential boundary layer 
that corresponds to only high-order error, i.e., boundary layers that only appear when measuring the 
approximation error in the $W^{1,\infty}$- or $H^1$-semi-norms, when defined.

In order to motivate our choice of boundary flux values when $r=0$, 
we observe that, for this special case, 
the DG gradient approximations $q^\pm_h$ are actually equivalent 
to the forward and backward difference quotients used in FD methods  
for interior simplexes 
when $\mathcal{T}_h$ is a Cartesian partition labelled with the natural ordering  
(see \cite{Feng_Lewis_Neilan13}).  
By extending the equivalence of the proposed LDG methods and the FD methods 
defined in \cite{Feng_Lewis12a,Lewis_Dissertation} to the boundary of the domain, 
we can derive the necessary boundary flux values for $u_h$ and $q_h^\pm$ on $\cE^B_h$.  
To this end, we will need to develop 
a methodology for extending the solution $u$ to the exterior of the domain $\Omega$.  
We now define a way to do such an extension that is consistent with the interpretation of the 
auxiliary variables and consistent with the FD strategy of introducing 
``ghost values" for a grid function, 
where the underlying grid will be defined by the midpoints of the Cartesian partition $\mathcal{T}_h$.  

We first describe the extension for the approximation function $u_h$.  
Given the Dirichlet boundary data for the viscosity solution $u$, it is natural to assume that the 
approximation function $u_h$ has a constant extension beyond each individual boundary face/edge.
Thus, we wish to define the exterior boundary fluxes using the Dirichlet boundary condition by setting 
$u ( x^E ) = g(x) $ for all $x \in \partial K \cap \partial \Omega$.  
However, a given boundary simplex may have multiple faces/edges in $\cE^B_h$.
Therefore, we introduce a ``ghost simplex" exterior to each individual face/edge in $\cE^B_h$, 
and we define the exterior value as $g_e$, where
\begin{equation} \label{DG_u_ext1}
	 \sum_{i=1}^d \bigl\langle g_e , n_e^{(i)} \bigr\rangle_{e} 
	 	= \sum_{i=1}^d \bigl\langle g , n_e^{(i)} \bigr\rangle_{e} \qquad\forall e\in \cE^B_h.
\end{equation} 
Then, we define 
\begin{equation} \label{DG_u_ext}
	u_h(x^E)\bigl|_{e} := g_e   \qquad\forall e\in \cE^B_h.
\end{equation}

Observe that, for $r = 0$, 
we only apply the Dirichlet boundary condition to the exterior function limits.
Furthermore, we define the exterior function limits to be edge-dependent.  
Since the function value is constant on each simplex $K$, we do not extend the Dirichlet boundary
condition to the interior of the domain by strongly enforcing \eqref{bc_ell}.  
Instead, we treat the value of $u_h$ on $K$ as an unknown whenever $K$ is a boundary simplex.
We use the edge-dependent definition to mimic the use of ghost values when $r=0$, 
which are introduced for each coordinate direction when using a FD methodology.  
When $\cT_h$ is a Cartesian partition, our methodology does in fact result in the introduction of a 
fixed exterior boundary flux value for each individual coordinate direction.  
The result of the methodology will be a more weighted approximation on a boundary simplex 
based upon the boundary condition along each boundary face/edge independently 
and on the PDE for the interior of the simplex.

Next, we describe how we assign boundary values for $q_h^\pm$  for $r = 0$.  
Since we do not have Neumman boundary data, we will have to enforce auxiliary boundary conditions.  
Assuming $\mathcal{T}_h$ is a Cartesian partition labelled with the natural ordering, throughout the interior of the domain there holds  
\begin{equation} \label{DG_q+-}
	q_i^- \bigl|_K = q_i^+ \bigl|_{K_i^-} , \qquad 
	q_i^+ \bigl|_K = q_i^- \big|_{K_i^+}
\end{equation}
for all $i=1,2,\ldots,d$ and all interior simplexes $K \in \cT_h$ due to the equivalence with FD forward 
and backward difference quotients, 
where $K_i^-$ denotes the neighboring simplex in the negative $i$-th Cartesian direction 
and $K_i^+$ denotes the neighboring simplex in the positive $i$-th Cartesian direction. 
Extending \eqref{DG_q+-} to the boundary yields
\begin{subequations} \label{DG_q+-bc}
\begin{align}
	q_i^- (x^E) & = q_i^+(x^I) , \qquad \text{if } n_{e}^{(i)} < 0 ,  \\
	q_i^+ (x^E) & = q_i^-(x^I) , \qquad \text{if } n_{e}^{(i)} > 0 
\end{align}
\end{subequations}
for $x \in e$, 
where both $q_i^+(x^I)$ and $q_i^-(x^I)$ are treated as unknowns.
We will assume such a relationship holds along the boundary for all triangulations.  
We also note that the relationship is arbitrary if $n_i^{(i)} = 0$.  

Observe that the above extension does not define exterior limits for 
$q_i^+$ if $n_{e}^{(i)} < 0$ or $q_i^-$ if $n_{e}^{(i)} > 0$.
In order to define the remaining exterior limit values, we impose the following auxiliary constraint equations:
\begin{subequations} \label{DG_bc2ab}
\begin{align} 
\sum_{i = 1}^d \Bigl\langle q_i^- ( x^I ) - q_i^- ( x^E ) , n_{e}^{(i)} \Bigr\rangle_{e} & = 0 
\qquad \forall e \in \cE^B_h,  \label{DG_bc2a} \\
\sum_{i = 1}^d \Bigl\langle q_i^+ ( x^I ) - q_i^+ ( x^E ) , n_{e}^{(i)} \Bigr\rangle_{e} & = 0 
\qquad \forall e \in \cE^B_h.  \label{DG_bc2b}
\end{align}
\end{subequations}
 
The above constraint equations are consistent with discretizing the higher order 
auxiliary constraint for all ghost-values of $q_h^\pm$: 
\begin{align*} 
\sum_{k = 1}^d \bigl( q_k^\pm \bigr)_{x_k} ( x ) & = 0 \qquad \forall x \in \Omega^c.
\end{align*}
The philosophy for such an auxiliary assumption can be found in  \cite{Feng_Neilan08}.
We note that the constraint equations \eqref{DG_bc2ab}
are also trivially satisfied when defining the exterior values for $r \geq 1$ 
due to our ``continuity" assumption.  
Assuming that $\cT_h$ is either a uniform Cartesian partition 
or a $d$-triangular partition where each simplex has at most one face/edge in $\cE^B_h$, 
we can see that all exterior limits on the boundary of the domain 
have now been expressed in terms of unknown interior limits that correspond to degrees of 
freedom for the discretization.  

We end this section by explicitly specifying the resulting exterior limit definitions for $q^+_h$ and $q^-_h$ 
when approximating a two-dimensional problem with piecewise constant basis functions.
The explicit definitions for one-dimensional problems can be found in \cite{Feng_Lewis12c}.
Let  $q_i^\pm := ( q_h^\pm )_i$.  Then, using the strategy given above, we have
\begin{alignat*}{2}
	q_1^+ ( x^E ) & = q_1^- ( x^I ) , &&\qquad 	q_2^+ ( x^E )  = q_2^- ( x^I ) , \\ 
	q_1^- ( x^E ) & = q_1^- ( x^I ) ,  &&\qquad	q_2^- ( x^E ) = q_2^- ( x^I ) ,
\end{alignat*}
if $n_1(x) < 0$ and $n_2(x) < 0$, 
\begin{alignat*}{2}
	q_1^+ ( x^E ) & = q_1^-( x^I ) , &&\qquad 	q_2^+ ( x^E )  = q_2^+ ( x^I ) 
		+ q_1^+ ( x^I ) - q_1^+ ( x^E ) , \\
	q_1^- ( x^E ) & = q_1^- ( x^I) 
		+ q_2^- ( x^I ) - q_2^- ( x^E ) , &&\qquad 	q_2^- ( x^E )  = q_2^+ ( x^I ) 
\end{alignat*}
if $n_1(x) < 0$ and $n_2(x) \geq 0$, 
\begin{alignat*}{2}
	q_1^- ( x^E) & = q_1^+ ( x^I ) , &&\qquad 	q_2^- ( x^E ) = q_2^- ( x^I ) + q_1^- ( x^I ) - q_1^- ( x^E ) , \\
	q_1^+ ( x^E ) & = q_1^+ ( x^I )  + q_2^+ ( x^I ) - q_2^+ ( x^E ) , &&\qquad 	q_2^+ ( x^E ) = q_2^-( x^I)  
\end{alignat*}
if $n_1(x) \geq 0$ and $n_2(x) < 0$, 
and 
\begin{alignat*}{2}
	q_1^- ( x^E ) & = q_1^+ ( x^I ) , &&\qquad 	q_2^-( x^E )  = q_2^+ ( x^I ) , \\ 
	q_1^+ ( x^E ) & = q_1^+( x^I ) , &&\qquad 	q_2^+ ( x^E )  = q_2^+( x^I ) 
\end{alignat*}
if $n_1(x) \geq 0$ and $n_2(x) \geq 0$ for all $x \in \partial \Omega \cap e$ for some $e \in \cE^B_h$.

\begin{remark} \ 
\begin{enumerate}[(a)]
\item 
When $r=0$, our approximation space consists of totally discontinuous 
piecewise constant functions.  
We have prescribed a way to assign all exterior boundary flux values for our approximation functions, 
and, by convention, we treat all interior boundary flux values as unknowns.
\item 
The above constraint equations occur naturally in the boundary edge terms 
for the bilinear form \eqref{a_bil} for each auxiliary variable.  
We use this observation to enforce our boundary conditions for $u_h$ and $q_h^\pm$ 
in the numerical tests found in section~\ref{sec-5}.     
\end{enumerate}
\end{remark}  

\subsection{The numerical viscosity and numerical moment} \label{sec-3_2b}

In this section, we take a closer look at the numerical viscosity and the numerical moment 
used in the definition of the Lax-Friedrichs-like numerical operator \eqref{LF1}.  
We divide the analysis into two cases, $r=0$ and $r \geq 1$.  
When $r = 0$, we will recover vanishing FD approximations 
of the Laplacian operator and the biharmonic operator.  
When $r \geq 1$, we will recover interior jump/stabilization terms.

First we consider the case $r = 0$ in the definition of $V^h$. Suppose that 
$\mathcal{T}_h$ is a uniform Cartesian partition labelled using the natural ordering.
Let $K$ be an interior simplex, $x_K$ denote its midpoint, and 
$\chi_K$ denote the characteristic function on $K$. Then, by \cite{Feng_Lewis_Neilan13}, 
we have
\begin{align*} 
	- \beta \cdot \bigl( q_h^+ - q_h^- , \chi_K \bigr)_{\cT_h} 
	& = - \sum_{i=1}^d \beta_i \bigl( \delta_{x_i,h_i}^+ u_h(x_K) - \delta_{x_i,h_i}^- u_h(x_K) \bigr) \\ 
	& = \sum_{i=1}^d \beta_i h_i \delta_{x_i, h_i}^2 u_h(x_K), 
\end{align*}
where $\delta_{x_i, h_i}^+$ denotes the forward difference quotient operator, 
$\delta_{x_i, h_i}^-$ denotes the backward difference quotient operator, 
and $\delta_{x_i, h_i}^2$ denotes the standard second order central difference quotient operator 
for approximating pure second derivatives.  
Also, by \cite{Feng_Lewis_Neilan13}, we have 
\begin{align*} 
	& \alpha : \bigl( P_{i j}^{+ +} - P_{i j}^{+ -} - P_{i j}^{- +} + P_{i j}^{- -} , \chi_K \bigr)_{\cT_h} \\
	& \qquad = \sum_{i,j=1}^d \alpha_{i j} 
		\bigl( \delta_{x_i, h_i}^+ \delta_{x_j,h_j}^+ u_h(x_K) 
		- \delta_{x_i, h_i}^+ \delta_{x_j,h_j}^- u_h(x_K) \\ 
	& \qquad \qquad \qquad 
		- \delta_{x_i, h_i}^- \delta_{x_j,h_j}^+ u_h(x_K) 
		+ \delta_{x_i, h_i}^- \delta_{x_j,h_j}^- u_h(x_K) \bigl) \\ 
	& \qquad = \sum_{i,j=1}^d \alpha_{i,j} h_i h_j \delta_{x_i, h_i}^2 \delta_{x_j, h_j}^2 u_h(x_K) .
\end{align*}
Thus, for $\beta = \vec{1}$ and $\alpha = \mathbf{1}$, we recover scaled approximations for the 
Laplace and biharmonic operator.  
Consequently, the Lax-Friedrichs-like numerical operator is a direct realization of the vanishing 
moment method (cf. \cite{Feng_Neilan08,Feng_Neilan11}) 
combined with the vanishing viscosity method from Hamilton-Jacobi equations (cf.  \cite{Crandall_Lions83}).  

A similar consequence of the relationship with FD when $r=0$ and $\mathcal{T}_h$ 
corresponds to a uniform Cartesian grid labelled using the natural ordering is that 
\[
	\bigl( P^{\pm \mp}_h \bigr)_{i i} = \frac{1}{h_i} \bigl( q_h^- - q_h^+ \bigr)_i \qquad\mbox{for } i=1,2,\cdots, d.
\]
Thus, if $\hF$ is defined by \eqref{LF1}, then $\hF$ may implicitly be monotone increasing 
with respect to $q_h^+$ and monotone decreasing with respect to $q_h^-$ 
for $\beta = \vec{0}$ as long as $h_i$ is sufficiently small and $\alpha_{i i} > 0$ 
for all $i = 1, 2, \ldots, d$.  
In other words, the numerical moment can implicitly enforce the g-monotonicity 
requirements for $q_h^\pm$.  
We exploit this observation in section~\ref{sec-5} by choosing $\beta = \vec{0}$ 
in our numerical tests.   
Heuristically, we expect the corresponding FD schemes to be limited to 1st order 
accuracy when the numerical viscosity is present (as with Lax-Friedrichs schemes for 
Hamilton-Jacobi equations), 
whereas the corresponding FD schemes may be capable of 2nd order accuracy 
when only the numerical moment is present.  
Such an observation is supported by the numerical tests found later in section~\ref{sec-5} 
as well as the numerical tests of the FD methods found in \cite{Lewis_Dissertation}.  

We now consider the case $r \geq 1$ in the definition of $V^h$.
Let $i \in \{ 1, 2, \ldots, d \}$. 
Observe that by the boundary conditions from section~\ref{sec-3_1}, we have 
\begin{align*}
	\bigl( q_i^+ - q_i^- , \phi \bigr)_{\cT_h}
	= a_i^+ \left( u_h , \phi \right) - a_i^- \left( u_h , \phi \right) 
	= \Big\langle \bigl[ u_h \bigl] , \bigl[ \phi \bigr] \, \big| n_{e}^{(i)} \big| \Big\rangle_{\cE^I_h}.
\end{align*}
Thus, 
\begin{equation} \label{DG_num_visc}
	- \beta \cdot \bigl( q_h^- - q_h^+ , \phi \bigr)_{\cT_h} 
	= \sum_{i=1}^d \beta_i \Big\langle \big[ u_h \big] , 
		\big[ \phi \big] \, \big| n_{e}^{(i)} \big| \Big\rangle_{\cE^I_h}. 
\end{equation}
Similarly, for $i,j \in \{ 1, 2, \ldots, d \}$, 
\begin{align*}
	& \bigl( P_{i,j}^{+ +} - P_{i,j}^{+ -} - P_{i,j}^{- +} + P_{i,j}^{- -} , \phi \bigr)_{\cT_h} \\ 
	& \qquad = a_j^+ \bigl( q_i^+ , \phi \bigr) - a_j^- \bigl( q_i^+ , \phi \bigr) - a_j^+ \bigl( q_i^- , \phi \bigr) 
		+ a_j^- \bigl( q_i^- , \phi \bigr) \\ 
	& \qquad = \Big\langle \big[ q_i^+ \big] , \bigl[ \phi \bigr] \, \big| n_{e}^{(j)} \big| \Big\rangle_{\cE^I_h} 
		- \Big\langle \big[ q_i^- \big] , \bigl[ \phi \bigr] \, \big| n_{e}^{(j)} \big| \Big\rangle_{\cE^I_h} .
\end{align*}
Thus, 
\begin{equation} \label{DG_num_mom}
	\alpha : \bigl( P_{i,j}^{+ +} - P_{i,j}^{+ -} - P_{i,j}^{- +} + P_{i,j}^{- -} , \phi \bigr)_{\cT_h}
	 = \sum_{i,j=1}^d \alpha_{i,j} \Big\langle \big[ q_i^+ - q_i^- \big] , 
	 	\bigl[ \phi \bigr] \, \big| n_{e}^{(j)} \big| \Big\rangle_{\cE^I_h}.
\end{equation}

From above, we can see that 
\begin{align}
	& \bil{a}_0 \bigl(u_h , q^-_h, q^+_h, P^{- -}_h, P^{- +}_h, P^{+ -}_h, P^{+ +}_h; \phi_h \bigr) \label{DG_fhat_rewritten} \\ 
	\nonumber
	& \qquad = \bigl( F \left( P_h , q_h , u_h , \cdot \right) , \phi_h \bigr)_{\cT_h} 
		+ \sum_{i=1}^d \beta_i \Big\langle \bigl[ u_h \bigr] , 
			\bigl[ \phi_h \bigr] \, \big| n_{e}^{(i)} \big| \Big\rangle_{\cE^I_h} \\ 
	\nonumber
		& \qquad \qquad 
		+ \sum_{i,j=1}^d \alpha_{i,j} \Big\langle \bigl[ q_i^+ - q_i^- \bigr] , 
			\bigl[ \phi_h \bigr] \, \big| n_e^{(j)} \big| \Big\rangle_{\cE^I_h} , 
\end{align}
where 
\[
	P_h = \frac{P_h^{+ -} + P_h^{- +}}{2} , \qquad q_h = \frac{q_h^+ + q_h^-}{2}, 
\]
and $q^+_h$, $q^-_h$ are both approximations for $\nabla u$.
Thus, adding a numerical moment and a numerical viscosity 
amounts to the addition of interior jump/stabilization terms to 
an $L^2$-projection of the fully nonlinear PDE operator into $V^h$.
We do note that the jump/stabilization terms that arise due to the numerical moment 
penalize the differences in $q_h^+$ and $q_h^-$.
Thus, the numerical moment is not analogous to a high order penalization term 
that penalizes jumps in a single approximation for $\nabla u$, as sometimes used in 
interior penalty methods.  
Instead, the numerical moment penalizes the difference in two optimal DG approximations for
 $\nabla u$ (cf. \cite{Feng_Lewis_Neilan13}).
We remark that this new jump term is the distinguishing characteristic of the proposed 
LDG methods since it was not possible to obtain an analogous result for the IPDG framework 
proposed in \cite{Feng_Lewis13}.

\subsection{Solvers} \label{sec-3_4}

We now discuss different strategies for solving the nonlinear system of equations 
that results from the proposed LDG discretization for the elliptic problem.  The underlying goal for the 
methodology presented in this paper is to discretize the fully nonlinear PDE problem 
in a way that removes much of the burden of approximating viscosity solutions from the design of 
the solver.  Thus, our primary focus is at the discretization level.  However, some of the properties 
of the methodology are more apparent from the solver perspective.

Most tests show that it is sufficient to simply use a Newton solver on the full system of 
equations  \eqref{pde_ell_weak} and \eqref{global_aux}.
Observe that only \eqref{pde_ell_weak} is nonlinear, the equation is purely algebraic, and $\hF$ is 
monotone in seven of its arguments.
The auxiliary equations \eqref{global_aux} are all linear.
The numerical operator presented in this paper is symmetric in both the 
mixed approximations $P_h^{- +}$ and $P_h^{+ -}$ and the non-mixed approximations $P_h^{- -}$ and 
$P_h^{+ +}$.  Thus, we can reduce the size of the system of equations by averaging the two 
pairs of auxiliary variables in the above formulation without changing the methodology.

Due to the size of the mixed formulation, we first present a splitting algorithm that 
provides an alternative to a straightforward Newton solver for the entire system of equations.
By using a splitting algorithm, the resulting algorithm will iteratively solve an entirely local, 
nonlinear equation 
that has strong monotonicity properties in the $d$ unknown arguments, 
and the solution of the equation can be mapped to an updated approximation for $u_h$.
Tests show that the solver is particularly useful for nonlinear problems 
that have a unique viscosity solution only defined in a restrictive function class.  
For instance, viscosity solutions of the Monge-Amp\'ere equation are unique in the class of convex functions.  
However, the proposed solver is not as efficient 
as the second solver we present that takes advantage of the above nonstandard discretization 
technique.  In order to improve the speed of the solver, fast Poisson solvers for the 
DWDG method (cf. \cite{Lewis_Neilan13}) need to be developed.

Our second solver strategy is a natural generalization of the FD 
methodology for numerical PDEs.
Constructing and applying the DG derivative operators requires 
sparse matrix multiplication and addition as well as inverting the local mass matrices.
Thus, all auxiliary equations in the mixed formulation can be solved for a given function $u_h$.
Substituting these operators directly into the numerical operator results in a single nonlinear 
variational problem for $u_h$ that can be solved iteratively.


\subsubsection{An inverse-Poisson fixed-point solver} \label{DG_DWDG_Solver}

We now describe the above mentioned splitting algorithm that takes into account 
the special structure of the nonlinear algebraic system that results from our nonstandard LDG 
discretization methods for elliptic PDEs and parabolic PDEs when using implicit time-stepping.
The algorithm is strongly based upon using a particular numerical moment.

\begin{algorithm}\ \label{DG_solver_alg}
\begin{enumerate}
\item Pick an initial guess for $u_h$. 
\item Form initial guesses for $q^+_h$, $q^-_h$, $P^{+ +}_h$, $P^{+ -}_h$, $P^{- +}_h$, and 
$P^{- -}_h$ using equations \eqref{global_aux}.  
\item Set
\begin{align*}
	G_i & := F \Bigl( \frac{P_h^{- +} + P_h^{+ -}}{2} , \frac{ q_h^- + q_h^+}{2} , u_h , x \Bigr)
		+ \gamma \bigl( P_h^{+ +} - P_h^{+ -} - P_h^{- +} + P_h^{- -} \bigr)_{i i} \\ 
		& \qquad - \beta_i \bigl( q_h^- - q_h^+ \bigr)_i 
\end{align*}
for a fixed constant $\gamma > 0$, and solve
\[
	\bigl( G_i , \varphi_{i} \bigr)_{\cT_h} = 0 \qquad \forall \varphi_{i} \in V^h
\]
for $ \frac12 \bigl( P_h^{- +} + P_h^{+ -} \bigr)_{i i}$ for all $i = 1, 2, \ldots, d$.
For sufficiently large $\gamma$ and a differentiable operator $F$, 
the above set of equations has a negative definite Jacobian.
\item Find $u_h$, $q_h^+$, and $q_h^-$ by solving the linear system of equations formed by 
\eqref{q_bil} and the trace of averaging \eqref{p_bil} for $\mu=-, \nu=+$ and $\mu=+, \nu=-$.
Observe that this is equivalent to solving Poisson's equation with source data given by the 
trace of $\frac12 \bigl( P_h^{- +} + P_h^{+ -} \bigr)$.
Alternatively, apply the DWDG method using 
the trace of $\frac12 \bigl( P_h^{- +} + P_h^{+ -}\bigr) $ as the source data to find $u_h$. 
\item Solve \eqref{p_bil} for $P_h^{+ +}$, $P_h^{+ -}$, $P_h^{- +}$, and $P_h^{- -}$.
If the alternative approach in step 4 was used, also solve \eqref{q_bil} for $q_h^+$ and $q_h^-$.  
\item Repeat Steps 3 - 5 until the change in $\frac12 \bigr( P_h^{- +} + P_h^{+ -}\bigr)$ is sufficiently small.
\end{enumerate}
\end{algorithm}

\noindent 
We now make a couple of comments about the proposed solver.

\begin{remark}\
\begin{enumerate}[(a)]
\item
The proposed algorithm is well-posed since it is based on the DWDG method which 
results in a symmetric positive definite discretization of Poisson's equation (cf. \cite{Lewis_Neilan13}).
\item
The nonlinear equation in Step 3 is entirely local with respect to the unknown variable.
\item
Clearly a fixed point for the solver corresponds to a discrete solution of the original 
PDE problem. In section~\ref{sec-5} and in \cite{Feng_Lewis12c}, 
we demonstrate that the above 
solver can be used to eliminate numerical artifacts that arise due to low-regularity PDE 
artifacts.  
Thus, the proposed solver is less dependent upon the initial guess.
The algorithm can also be used to form a preconditioned initial guess for other nonlinear solvers 
that may be faster but require a ``better" initial guess.  
\end{enumerate}
\end{remark}


\subsubsection{A direct approach for a reduced system} \label{DG_direct_Solver}

In this section, we propose a solver technique that is analogous to the approach used in FD methods.  
Observe that if 
$\bigl(u_h , q^+_h, q^-_h, P^{+ +}_h, P^{+ -}_h, P^{- +}_h, P^{- -}_h \bigr)$
is a solution to \eqref{pde_ell_weak} and \eqref{global_aux}, 
then there exists linear operators $\nabla_h^\pm$ and $D_h^{\mu \nu}$ 
such that $q^\pm_h = \nabla_h^\pm u_h$ and 
$P_h^{\mu \nu} = D_h^{\mu \nu} u_h$ for all $\mu, \nu \in \{+, - \}$, 
where the linear operators are locally defined by \eqref{DG_local_q} and \eqref{DG_local_p}.  

Using these numerical derivative operators, the second solver is given by: 
 \begin{algorithm}\ \label{DG_fd_solver}
 \begin{enumerate}
 \item Given $\cT_h$ and $V^h$, compute the operators $\nabla_h^{\pm}$ and $D_h^{\mu \nu}$. \\
 \item Solve for $u_h \in V^h$ the single nonlinear equation
\begin{align*}
	\Bigl( \hF \bigl( D_h^{+ +} u_h , D_h^{+ -} u_h , D_h^{- +} u_h , D_h^{- -} u_h , 
		\nabla_h^{+} u_h , \nabla_h^{-} u_h , u_h , \cdot \bigl) , \varphi_h \Bigr)_{\cT_h}=0 \quad
		\forall \varphi_h \in V^h.
\end{align*}
 \end{enumerate}
 \end{algorithm}
 
 We note that a reduced formulation can also be used where we simply create 
 the following new differential operators:
 \[
	 \overline{D}_h^2 := \frac{ D_h^{- -} + D_h^{+ +} }{2}, \qquad 
	 \widetilde{D}_h^2 := \frac{ D_h^{- +} + D_h^{+ -} }{2}, \qquad 
	 \nabla_h := \frac{\nabla^+_h + \nabla^-_h}{2} . 
 \] 
 The Lax-Friedrichs-like numerical operator can be witten as 
 \begin{align} \label{reduced_direct}
 	& \hF \bigl( \overline{D}_h^2 u_h , \widetilde{D}_h^2 u_h , \nabla_h^{+} u_h , \nabla_h^{-} u_h , u_h , x \bigr) \\ 
	& \quad \nonumber
		= F \bigl( \widetilde{D}_h^2 u_h , \nabla_h u_h , u_h, x \bigr)  
		+ 2 \alpha : \big( \overline{D}_h^2 u_h - \widetilde{D}_h^2 u_h \big) 
		- 2 \beta \cdot \big( \nabla_h^+ u_h - \nabla_h^- u_h \big).
 \end{align}
 
\noindent 
For all of the tests below where a Newton solver is used for the full system of 
equations in the mixed formulation, analogous results were obtained using 
Algorithm~\ref{DG_fd_solver} with the reduced numerical operators.
As expected, for two-dimensional problems we observed significant speed-up in the 
performance of the solver.

\begin{remark}
The methodology of Algorithm~\ref{DG_fd_solver} follows directly from the FD methodology 
where derivatives in a PDE are simply replaced by numerical derivatives of 
the approximation for the solution $u$ to form the discretization of the PDE problem.  
For nonlinear problems, we replace the nonlinear PDE operator by a numerical operator.
In our LDG setting, we use the LDG methodology to define the various numerical derivatives. 
\end{remark}

\section{An extension for parabolic problems} \label{sec-4}

We now develop fully discrete methods for approximating the parabolic equation 
\eqref{pde} complemented by the following boundary condition and initial condition: 
\begin{subequations} 
\begin{alignat}{2}
	u(x,t) & = g(x), \qquad && (x,t) \in \Omega_T := \Omega \times (0,T] , 
		\label{bc_time} \\ 
	u(x,0) & = u_0(x) , \qquad && x \in \Omega 
		\label{ic_time}
\end{alignat}
\end{subequations}
using an LDG spatial-discretization paired with the method of lines approach for the time discretization.
Taking advantage of the elliptic formulation in section~\ref{sec-3}, we will propose 
the following implicit and explicit time-discretizations:
forward Euler, backward Euler, trapezoidal, and Runge-Kutta (RK).
The time-discretization used in application should be selected according to the potential optimal order 
$r+1$ of the LDG spatial-discretization for sufficiently regular viscosity solutions.

We first present the semi-discrete discretization of the (fully) nonlinear equation \eqref{pde}
by discretizing the spatial dimension.
Replacing the PDE operator $F$ with a numerical operator $\hF$ in \eqref{pde},  
applying a spatial discretization using the above LDG framework for elliptic equations, 
and using the $L^2$-projection operator $\cP_h : L^2(\cT_h) \to V^h$ defined by
\begin{equation} \label{DG_proj}
	\bigl( \cP_h v , \phi_h \bigr)_{\cT_h} = \big( v , \phi_h \bigr)_{\cT_h} \qquad \forall \phi_h \in V^h 
\end{equation}
for all $v \in L^2(\cT_h)$, 
we have the 
following semi-discrete equation
\begin{equation}\label{DG_semi-discrete}
	{(u_h)}_t = - \cP_h \Bigl( \hF \bigl( P_h^{+ +}, P_h^{+ -}, P_h^{- +}, P_h^{- -}, q_h^{+}, q_h^{-}, u_h, x, t \bigr) \Bigr) , 
\end{equation}
where, given $u_h$ at time $t$, corresponding values for $q_h^{\pm}$ and $P_h^{\mu \nu}$, 
$\mu, \nu \in \{+, -\}$,  can be found by solving 
the local equations \eqref{DG_local_q} and \eqref{DG_local_p}. 

Our full-discretization of the initial-boundary value problem 
\eqref{pde}, \eqref{bc_time}, and \eqref{ic_time} is defined by applying an ODE solver
to the semi-discrete (variational) form given in \eqref{DG_semi-discrete}.
To partition the time domain, we fix an integer $M>0$ and 
let $\Delta t = \frac{T}{M}$.  
Then, we define $t_k := k \, \Delta t$ for a real number $k$ with $0 \leq k \leq M$.
Notationally, $u_h^k \in V^h$ and $q_h^{\pm, k} \in (V^h)^d$ 
will be an approximation for $u(\cdot, t_k)$ and 
$\nabla u(\cdot, t_k)$, respectively, for all $0 \leq k \leq M$.
For both implicit and explicit schemes, we define the initial value, $u_h^0$, by
\begin{equation} \label{DG_proj_ic}
	u_h^0 = \cP_h u_0 .
\end{equation}

To simplify the appearance of the methods and to make them more transparent for use with
a given ODE solver, we use a subscript $k$ to denote 
the fact that the boundary values are being naturally enforced in \eqref{DG_local_q} and \eqref{DG_local_p} 
using the boundary condition \eqref{bc_time} evaluated at time $t_k$, $0 \leq k \leq M$.  
Thus, 
\begin{align} \label{DG_q_time}
\Bigl( (q_{h,k}^\pm)_i , \phi^\pm_i \Bigr)_{\cT_h} 
	& =  \Big\langle T_i^\pm (u_{h,k}) , [ \phi^\pm_i ] \, n_e^{(i)} \Big\rangle_{\cE^I_h} 
	+ \Big\langle T_i^\pm (u_{h,k}) , \phi^\pm_i(x^I) \, n_i \Big\rangle_{\cE^B_h} \\
	\nonumber & \qquad 
	- \bigl( u_{h,k} , ( \phi^\pm_i )_{x_i} \bigr)_{\mathcal{T}_h} \qquad \forall \phi^\pm_i \in V^h
\end{align}
for $i = 1, 2, \ldots, d$, where we evaluate the boundary flux values 
using the convention 
\[
	\sum_{i=1}^d \bigl\langle u_{h,k} , \varphi_h \, n_i \bigr\rangle_{\cE^B_h} 
		= \sum_{i=1}^d \bigl\langle g(\cdot, t_k) , \varphi_h \, n_i \bigr\rangle_{\cE^B_h} 
		\qquad \forall \varphi_h \in V^h 
\]
when $r \geq 1$   and 
\[ 
	\sum_{i=1}^d \bigl\langle u_{h,k}(x^E) , n_e^{(i)} \bigr\rangle_{e} 
		= \sum_{i=1}^d \bigl\langle g(\cdot, t_k) , n_e^{(i)} \bigr\rangle_{e}   
\] 
when $r = 0$. Similarly, 
\begin{align} \label{DG_p_time}
\Bigl( \bigl(P^{\mu \nu}_{h,k} \bigr)_{i j} ,  \psi^{\mu \nu}_{i j}  \Bigr)_{\cT_h} 
& = \Big\langle T_j^\nu \bigl( (q_{h,k}^\mu)_i \bigr) , \bigl[\psi^{\mu \nu}_{i j} \bigr] \, n_e^{(j)} \Big\rangle_{\cE^I_h} 
	+ \Big\langle T_j^\nu  \bigl(q_{h,k}^\mu)_i \bigr) , \psi^{\mu \nu}_{i j} (x^I) \, n_j \Big\rangle_{\cE^B_h} \\ 
	\nonumber & \qquad - \bigl( \bigl(q_{h,k}^\mu \bigr)_i , 	(\psi^{\mu \nu}_{i j})_{x_j} \bigr)_{\mathcal{T}_h} 
\qquad \forall \psi^{\mu \nu}_{i j} \in V^h
\end{align}
for $i, j \in \{ 1, 2, \ldots, d \}$, $\mu, \nu \in \{ +, - \}$, 
where we assume $\bigl( q_{h,k}^\pm (x^E) \bigr)_i = \bigl( q_{h,k}^\pm (x) \bigr)_i$ when $ r \geq 1$ or 
\begin{align*}
	\sum_{i = 1}^d \Big\langle \bigl( q_{h,k}^\pm (x^I) \bigr)_i - \bigl( q_{h,k}^\pm (x^E) \bigr)_i, n_e^{(i)} \Big\rangle_{e} & = 0  
\end{align*}
and
\begin{align*}
	\bigl( q_{h,k}^- (x^E) \bigr)_i & = \bigl( q_{h,k}^+ (x^I) \bigr)_i , \qquad \text{if } n_e^{(i)} < 0 ,  \\
  \bigl( q_{h,k}^+ (x^E) \bigr)_i & =\bigl( q_{h,k}^- (x^I) \bigr)_i , \qquad \text{if } n_e^{(i)} > 0 
\end{align*}
for all $e \in \cE^B_h$, 
using \eqref{DG_q+-bc} and \eqref{DG_bc2ab}, when $r=0$.
Note, for $k=0$, we replace $g(\cdot, t_k)$ with $u_0(\cdot)$ in the above constraint equations 
if $u_0$ has an $L^2$ trace.  
Otherwise, we replace $g(\cdot, t_k)$ with the trace of $\mathcal{P}_h u_0$.  

We also simplify the presentation of the fully-discrete methods by introducing the operator notation 
\begin{equation}\label{DG_fhat_short}
\hF^k [v]  
:= \hF \left( D^{+ +}_{h,k} v , D^{+ -}_{h,k} v , D^{- +}_{h,k} v , D^{- -}_{h,k} v, \nabla^{+}_{h,k} v , \nabla^{-}_{h,k} v , 
	v , x, k \, \Delta t \right) 
\end{equation}
for all $v \in V^h$, 
where we are introducing linear operators $\nabla_{h,k}^\pm$ and $D_{h,k}^{\mu \nu}$ such that 
$q^\pm_{h,k} = \nabla_{h,k}^\pm u_h$ and 
$P_{h,k}^{\mu \nu} = D_{h,k}^{\mu \nu} u_h$ for all $\mu, \nu \in \{+, - \}$, 
where the linear operators are locally defined by replacing 
$u_{h,k}$ with an arbitrary function $v_h \in V^h$ in \eqref{DG_q_time} and \eqref{DG_p_time}.  
Then, the semi-discrete equation can be rewritten compactly as
\begin{equation}\label{DG_semi}
\bigl( u_h \bigr)_t (x, t_k) = - \cP_h \hF^k \bigl[ u_h(x, t_k) \bigr] \qquad \forall \, 0 \leq k \leq M,  x \in \Omega.
\end{equation}

Lastly, we define a modified projection operator 
$\cP_{h,k} : L^2(\cT_h) \to V^h$
that will be used to enforce the boundary conditions for explicit methods
using a penalty technique due to Nitsche in \cite{Nitsche70}.
Thus, we define $\cP_{h,k}$ by 
\begin{align} \label{DG_proj_bc}
	& \bigl( \cP_{h,k} v , \varphi_h \bigr)_{\cT_h} 
		+ \delta \sum_{i=1}^d \Big\langle \cP_{h,k} v , \varphi_h \, n_i \Big\rangle_{\cE^B_h} \\ 
	& \qquad 
	= \big( v , \varphi_h \bigr)_{\cT_h} 
		+ \delta \sum_{i=1}^d \Big\langle g(\cdot, t_k), \varphi_h \, n_i \Big\rangle_{\cE^B_h}  
		\qquad \forall \varphi_h \in V^h 
	\nonumber
\end{align}
for all $v \in L^2(\cT_h)$, 
where $\delta$ is a nonnegative penalty constant and $0 \leq k \leq M$.
We note that, for $\delta = 0$, $\cP_{h,k} = \cP_h$, yielding the broken $L^2$-projection operator.

Using the above conventions, we can define fully discrete methods for approximating problem 
\eqref{pde}, \eqref{bc_time}, and \eqref{ic_time}
based on approximating \eqref{DG_semi} using the forward Euler method, backward Euler method, 
or the trapezoidal method. Thus, we have respectively
\begin{equation} \label{DG_feuler}
	u_h^{n+1} = \cP_{h,n+1} \left( u_h^n - \Delta t \, \hF^{n} \left[ u_h^n \right] \right), 
\end{equation}
\begin{equation} \label{DG_beuler}
	u_h^{n+1} + \Delta t \, \cP_{h} \, \hF^{n+1} \left[ u_h^{n+1} \right] = u_h^n , 
\end{equation}
and
\begin{equation} \label{DG_trap}
	u_h^{n+1} + \frac{\Delta t}{2} \, \cP_{h} \, \hF^{n+1} \left[ u_h^{n+1} \right] = 
		u_h^n - \frac{\Delta t}{2} \, \cP_{h} \, \hF^{n} \left[ u_h^{n} \right]  
\end{equation}
for $n = 0, 1, \ldots, M-1$, 
where $u_h^0 := \cP_h u_0$ and, for \eqref{DG_beuler} and \eqref{DG_trap}, 
we also have, by \eqref{DG_fhat_short}, the implied auxiliary linear equations 
\begin{alignat*}{2}
q_h^{\mu , n } & = \nabla_{h,n}^\mu u_h^n  
	&& \qquad \forall \mu \in \{ +, - \} , \\
P_h^{\mu \nu , n} & = D_{h,n}^{\mu \nu} u_h^n  
	&& \qquad \forall \mu, \nu \in \{ +, - \}.
\end{alignat*}

\begin{remark}
Using an implicit method, such as the backward Euler and the trapezoidal method, results in approximating 
a fully nonlinear elliptic PDE at each time step using the LDG methods for elliptic PDEs formulated 
in section~\ref{sec-3}. Due to the time integration, the nonlinear solver has a natural initial guess for each time-step given 
by the approximation at the previous time step.
\end{remark}
 
Finally, we formulate the Runge-Kutta (RK) methods for approximating \eqref{DG_semi}.
 Let $s$ be a positive integer, $A \in \mathbf{R}^{s \times s}$, and $b,c \in \mathbf{R}^s$ such that
 \[
 	\sum_{\ell = 1}^s a_{k,\ell} = c_k
 \]
 for each $k = 1, 2, \ldots, s$.
 Then, a generic $s$-stage RK method for approximating \eqref{DG_semi} is defined by 
 \begin{equation} \label{DG_rk}
 	u_h^{n+1} = \cP_{h,n+1} \Bigl( u_h^n - \Delta t \sum_{\ell = 1}^s b_\ell \hF^{n+c_\ell} [ \xi_h^{n,\ell} ] \Bigr), \quad n = 0, 1, \ldots, N-1,
 \end{equation}
 where
 \[
 	\xi_h^{n,\ell} = \cP_{h, n+c_k} \Bigl( u_h^n - \Delta t \sum_{k = 1}^s a_{k,\ell} \hF^{n+c_k} [ \xi_h^{n,k} ] \Bigr),
 	\quad n = 0, 1, \ldots, N-1,
 \]
and  $u_h^0 = \cP_h u_0$.  We note that \eqref{DG_rk} corresponds to an explicit method when $A$ is strictly lower diagonal and 
 an implicit method otherwise.
 
 \begin{remark} 
$\xi_h^{n,\ell}$ in \eqref{DG_rk} can be viewed as an approximation for $u_h^{n+c_\ell}$.
 Since the boundary condition at $t_{n+1}$ is enforced by $\hF^{n+1}$, we can set $\delta = 0$ 
 in \eqref{DG_proj_bc} if $c_s = 1$.
 \end{remark}

\section{Numerical experiments} \label{sec-5}

In this section, we present a series of numerical tests to demonstrate 
the utility of the proposed LDG methods for fully nonlinear PDE problems of
type \eqref{pde_ell} and \eqref{pde} with two spatial dimensions.  
For elliptic problems, both Monge-Amp\`ere and Hamilton-Jacobi-Bellman 
types of equations will be tested. 
We also perform a test using the (semi-linear) 
infinite-Laplacian equation with a known low-regularity solution.  
The tests use spatial meshes composed of uniform rectangles.
To solve the resulting nonlinear algebraic systems, we use either 
the Matlab built-in nonlinear solver {\em fsolve} or Algorithm~\ref{DG_solver_alg}, where 
{\em fsolve} is used to perform Step 3 of Algorithm~\ref{DG_solver_alg}.
For the elliptic problems, we choose the initial guess as the zero function.
For the parabolic test problem, 
we choose the initial guess as the approximation formed at the previous time step 
and use the backward Euler method.
We also choose the approximation at time $t=0$ to be given by the $L^2$-projection 
of the initial condition into $V^h$.

For our numerical tests, errors will be measured in the $L^\infty$ norm 
and the $L^2$ norm. 
All recorded data corresponds to tests without a numerical viscosity, i.e., $\beta = \vec{0}$. 
Similar results hold when the numerical viscosity is present.   
For elliptic problems and parabolic problems where the error is not dominated by the time 
discretization, 
the test problems in \cite{Feng_Lewis12c} indicate 
the spatial errors are of order $\mathcal{O} (h^{s})$ for most problems, 
where $s = \min \{ r+1, k \}$ for the viscosity solution $u \in H^k(\Omega)$.  
In this paper, the computed convergence rates are a little more sporadic.  
On average, the schemes appear to exhibit an optimal rate of convergence in both norms.    
We note that the actual convergence rates have not yet been analyzed, 
and they may also depend on the regularity of the differential operator $F$ 
and the severity of its nonlinearity  
in addition to the regularity of the viscosity solution $u$.  


\begin{example} \label{DG_2D_test1}
Consider the Monge-Amp\`{e}re problem
\begin{align*}
- \textthm{det } D^2 u = -u_{x x} \, u_{y y} + u_{x y} \, u_{y x} &= f  \qquad \textthm{in } \Omega , \\
u &= g \qquad \textthm{on } \partial \Omega ,   
\end{align*}
where $f = -(1+x^2+y^2)e^{x^2+y^2}$, $\Omega = (0,1) \times (0,1)$, and $g$ is chosen such that
the viscosity solution is given by $u(x,y) = e^{\frac{x^2+y^2}{2}}$. 
\end{example}

\begin{figure}[htb]
\centering
\includegraphics[scale=0.35]{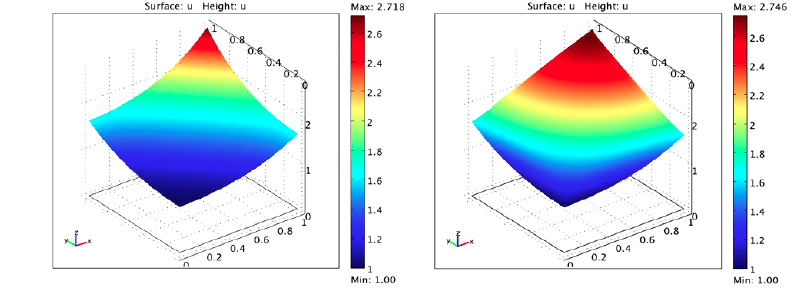}
\caption{
The two possible solutions for Example~\ref{DG_2D_test1}, as computed in \cite{Feng_Neilan08}.
The left plot corresponds to the viscosity solution while the right plot corresponds to the viscosity solution of 
$F[u] = \text{det }D^2 u$.
} 
\label{DG_macinf_soln}
\end{figure}

\begin{figure}[htb]
\centering
\includegraphics[scale=0.35]{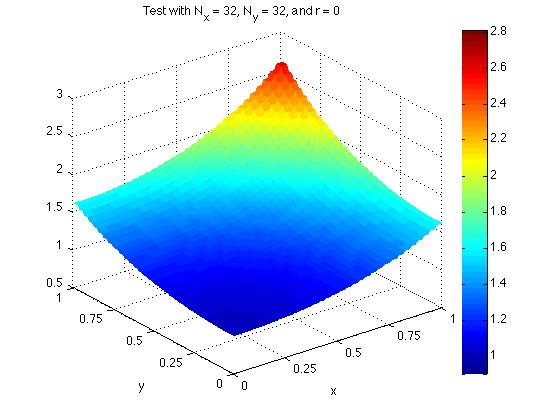}
\caption{
Computed solution for Example~\ref{DG_2D_test1} using 
$r = 0$, $\alpha = 24 I$, $h$ = 4.419e-02, and 
{\em fsolve} with initial guess $u_h^{(0)} = 0$. 
} 
\label{DG_ma_cinfty_r0_plot}
\end{figure}

\begin{figure}[htb]
\centering
\includegraphics[scale=0.35]{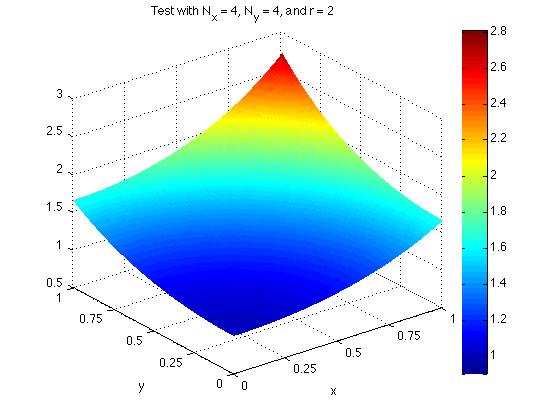}
\caption{
Computed solution for Example~\ref{DG_2D_test1} using 
$r = 2$, $\alpha = 24 I$, $h$ = 3.536e-01, and 
{\em fsolve} with initial guess $u_h^{(0)} = 0$. 
} 
\label{DG_ma_cinfty_r2_plot}
\end{figure}

Notice that the problem has two possible solutions as represented in Figure~\ref{DG_macinf_soln}.
Also, this problem is degenerate for the class of functions that are both concave and convex.
Results for approximating with $r=0,1,2$ can be found in Tables~\ref{DG_ma_cinfty_r0}, 
\ref{DG_ma_cinfty_r1}, and \ref{DG_ma_cinfty_r2}, respectively,
where we observe optimal convergence rates.
Plots for some of the various approximations can be found in 
Figures~\ref{DG_ma_cinfty_r0_plot} and \ref{DG_ma_cinfty_r2_plot}.

\begin{table}[htb]
\caption{
Rates of convergence for Example~\ref{DG_2D_test1} using 
$r = 0$, $\alpha = 24 I$, and 
{\em fsolve} with initial guess $u_h^{(0)} = 0$. 
}
\begin{center}
\begin{tabular}{| c | c | c | c | c |}
         \hline
    $h$ & $L^\infty$ norm & order & $L^2$ norm & order \\
         \hline
    1.41e-01 & 3.73e-01 &   & 8.31e-02 &   \\
         \hline
    8.84e-02 & 2.42e-01 &  0.92 & 5.10e-02 &  1.04   \\
         \hline
    5.89e-02 & 1.64e-01 &  0.95 & 3.31e-02 &  1.06   \\
         \hline
    4.42e-02 & 1.24e-01 &  0.97 & 2.44e-02 &  1.07   \\
         \hline
\end{tabular}
\end{center}
\label{DG_ma_cinfty_r0}
\end{table}

\begin{table}[htb]
\caption{
Rates of convergence for Example~\ref{DG_2D_test1} using 
$r = 1$, $\alpha = 24 I$, and 
{\em fsolve} with initial guess $u_h^{(0)} = 0$. 
}
\begin{center}
\begin{tabular}{| c | c | c | c | c |}
         \hline
    $h$ & $L^\infty$ norm & order & $L^2$ norm & order \\
         \hline
    1.41e-01 & 2.47e-02 &   & 1.73e-03 &   \\
         \hline
    1.18e-01 & 1.36e-02 &  3.25 & 1.61e-03 &  0.39   \\
         \hline
    1.01e-01 & 1.03e-02 &  1.81 & 1.12e-03 &  2.31   \\
         \hline
    7.86e-02 & 8.04e-03 &  0.99 & 5.82e-04 &  2.62   \\
         \hline
\end{tabular}
\end{center}
\label{DG_ma_cinfty_r1}
\end{table}

\begin{table}[htb] 
\caption{
Rates of convergence for Example~\ref{DG_2D_test1} using 
$r = 2$, $\alpha = 24 I$, and 
{\em fsolve} with initial guess $u_h^{(0)} = 0$. 
}
\begin{center}
\begin{tabular}{| c | c | c | c | c |}
         \hline
    $h$ & $L^\infty$ norm & order & $L^2$ norm & order \\
         \hline
    7.07e-01 & 6.39e-02 &   & 4.45e-03 &   \\
         \hline
    4.71e-01 & 2.32e-02 &  2.50 & 1.30e-03 &  3.03   \\
         \hline
    3.54e-01 & 1.09e-02 &  2.63 & 5.45e-04 &  3.02   \\
         \hline
\end{tabular}
\end{center}
\label{DG_ma_cinfty_r2}
\end{table}

\begin{figure}[htb]
\centering
\includegraphics[scale=0.35]{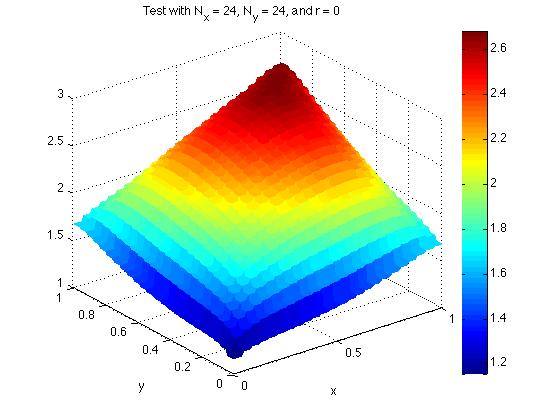}
\caption{
Computed solution for Example~\ref{DG_2D_test1} using 
$r=0$, $\alpha = -12 \mathbf{1}$, $h$ = 5.893e-02, and 
{\em fsolve} with initial guess $u_h^{(0)} = 0$.
}
\label{DG_mcinf_neg_moment}
\end{figure}

We now demonstrate that the numerical moment assists with resolving the issue of numerical artifacts 
and uniqueness only in a restrictive function class.  
We approximate Example~\ref{DG_2D_test1} using the numerical moment with $\alpha = -12 \mathbf{1}$, 
$N_x = N_y = 24$, $r=0$, and initial guess given by the zero function.  
The result is recorded in Figure~\ref{DG_mcinf_neg_moment}.
Thus, we can see that for a negative semi-definite choice for $\alpha$, we recover an 
approximation for the non-convex solution of the Monge-Amp\`{e}re problem represented in 
Figure~\ref{DG_macinf_soln}.  

\begin{example}  \label{DG_2D_test2}
Consider the Monge-Amp\`{e}re problem
\begin{align*}
- \textthm{det } D^2 u = -u_{x x} \, u_{y y} + u_{x y} \, u_{y x} &= 0 \qquad \textthm{in } \Omega , \\
u &= g \qquad \textthm{on } \partial \Omega, 
\end{align*}
where 
$\Omega = (-1,1) \times (-1,1)$ and $g$ is chosen such that
the viscosity solution is given by $u(x,y) = |x| \in H^1(\Omega)$.
\end{example}

\begin{figure}[htb]
\centering
\subfloat[$h_x$ = 6.667e-02.]{
\includegraphics[scale=0.35]{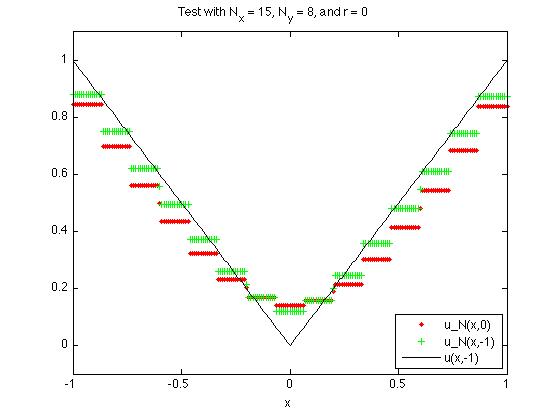}
}
\subfloat[$h_x$ = 1.818e-02.]{
\includegraphics[scale=0.16]{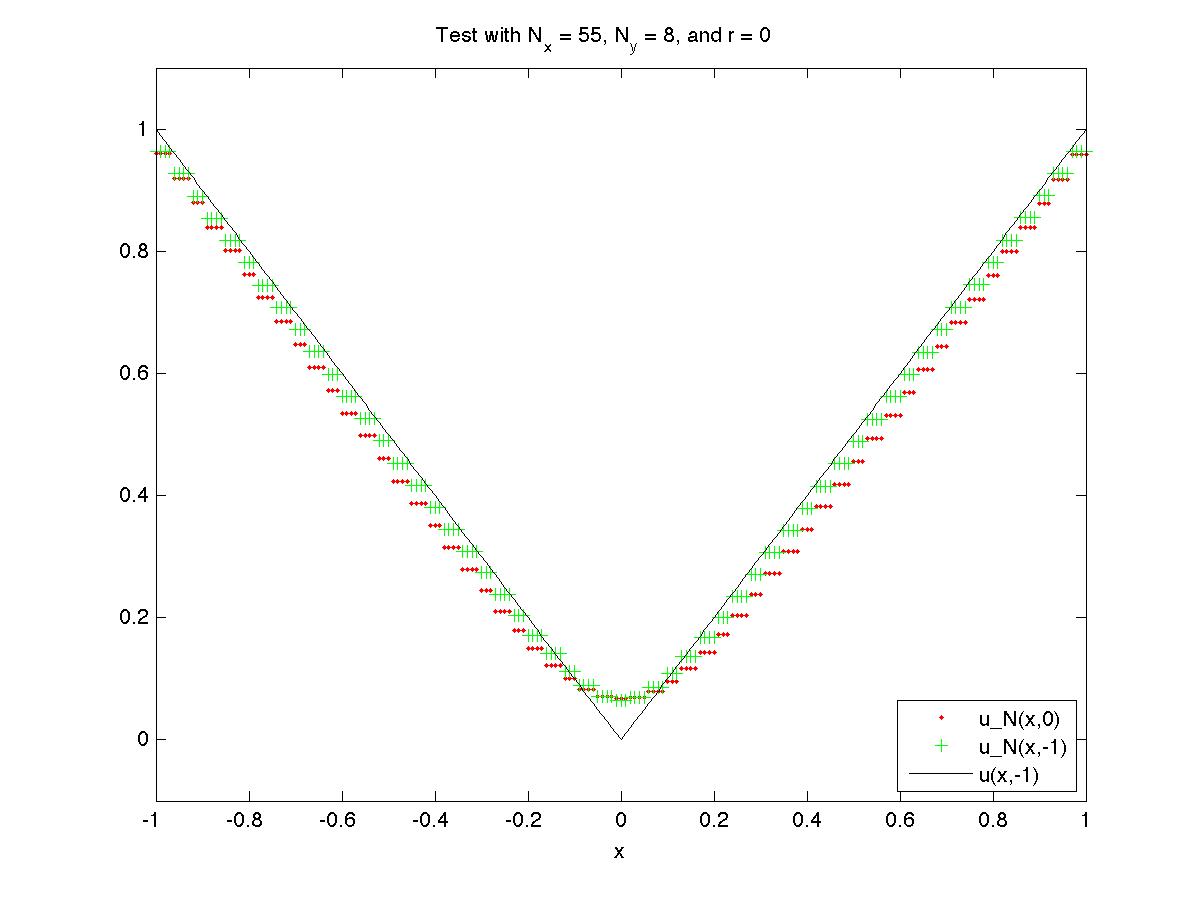}
} 
\caption{Computed solutions for Example~\ref{DG_2D_test2} using $r=0$, $\alpha = I$, $h_y$ = 1.250e-01, and 
{\em fsolve} with initial guess $u_h^{(0)} = 0$.}
\label{DG_ma_c0_r0_plot}
\end{figure}

\begin{figure}[htb]
\centering
\subfloat[$h_x$ = 4.167e-02 and $h_y$ = 1.667e-01.]{
\includegraphics[scale=0.35]{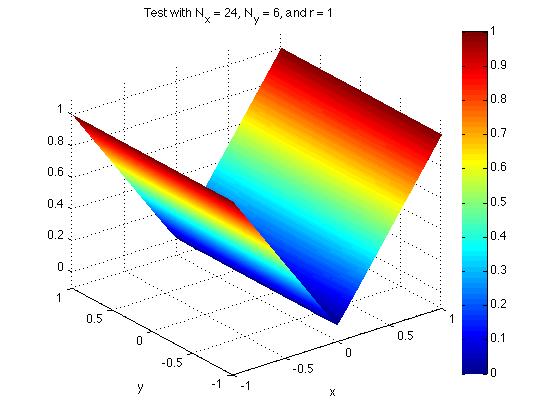}
}
\subfloat[$h_x$ = 4.167e-02 and $h_y$ = 1.667e-01.]{
\includegraphics[scale=0.35]{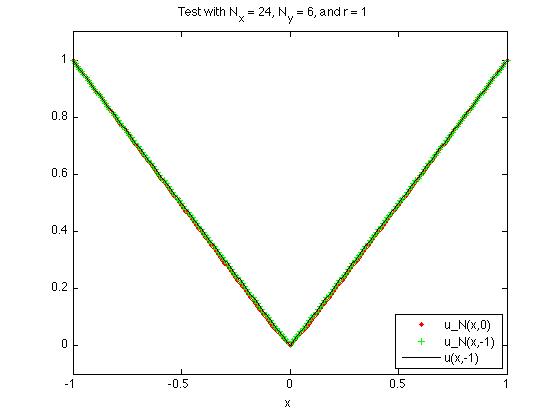}
}  \\ 
\subfloat[$h_x$ = 2.000e-01 and $h_y$ = 2.000e-01.]{
\includegraphics[scale=0.35]{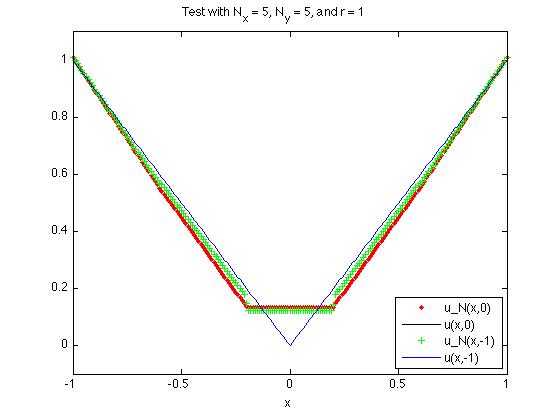}
}  
\caption{
Computed solution for Example~\ref{DG_2D_test2} using 
$r = 1$, $\alpha = I$, and 
Algorithm~\ref{DG_solver_alg} with initial guess $u_h^{(0)} = 0$.
Note that the top plots correspond to $x=0$ an edge and the bottom plot does not.}
\label{DG_ma_c0_r1_plot}
\end{figure}

Observe that the PDE is actually degenerate when acting on the solution $u$.
Furthermore, due to the low regularity of $u$, 
we expect the rate of convergence to be bound by one.
Using both piecewise constant and piecewise linear basis functions, we can see that the 
rate of convergence is bound by the theoretical bound in 
Table~\ref{DG_ma_c0_r0} and Table~\ref{DG_ma_c0_r1}.  
Plots for some of the approximations can be found in Figure~\ref{DG_ma_c0_r0_plot} for $r=0$ 
and Figure~\ref{DG_ma_c0_r1_plot} for $r=1$.
We remark that for $r=0$, all three solver approaches discussed in 
section~\ref{sec-3_4} gave analogous results.
However, for $r =1$, the direct formulation appears to have small residual wells that can trap the solver.
Thus, for this test, the non-Newton solver given by Algorithm~\ref{DG_solver_alg} 
appears to be better suited.

\begin{table}[htb]
\caption{
Rates of convergence for Example~\ref{DG_2D_test2} using 
$r=0$, $\alpha = I$, and 
{\em fsolve} with initial guess $u_h^{(0)} = 0$.
}
\begin{center}
\begin{tabular}{| c | c | c | c | c |}
         \hline
    $h_x$ & $L^\infty$ norm & order & $L^2$ norm & order \\
         \hline
    1.33e-01 & 1.87e-01 &   & 1.70e-01 &   \\
         \hline
    8.00e-02 & 1.30e-01 &  0.71 & 1.22e-01 &  0.65   \\
         \hline
    5.71e-02 & 1.02e-01 &  0.72 & 9.77e-02 &  0.66   \\
         \hline
    4.44e-02 & 8.51e-02 &  0.74 & 8.23e-02 &  0.68   \\
         \hline
    3.64e-02 & 7.33e-02 &  0.74 & 7.16e-02 &  0.69   \\
         \hline
\end{tabular}
\end{center}
\label{DG_ma_c0_r0}
\end{table}

\begin{table}[htb]
\caption{
Rates of convergence for Example~\ref{DG_2D_test2} using 
$r=1$, $\alpha = I$, $h_y = 1/3$ fixed, and 
Algorithm~\ref{DG_solver_alg} with initial guess $u_h^{(0)} = 0$.
}
\begin{center}
\begin{tabular}{| c | c | c | c | c |}
         \hline
    $h_x$ & $L^\infty$ norm & order & $L^2$ norm & order \\
         \hline
    2.50e-01 & 3.86e-02 &   & 3.42e-02 &   \\
         \hline
    1.25e-01 & 2.08e-02 &  0.89 & 1.85e-02 &  0.88   \\
         \hline
    8.33e-02 & 1.38e-02 &  1.02 & 1.24e-02 &  0.99   \\
         \hline
\end{tabular}
\end{center}
\label{DG_ma_c0_r1}
\end{table}

Another benefit of the numerical moment is that it can help regularize a problem that may
not be well-conditioned for a Newton solver due to a singular or poorly scaled Jacobian.  
Note that $\frac{\partial F}{\partial D^2 u} = 0$ almost everywhere in $\Omega$ for the
viscosity solution $u$ due to the fact that $D^2 u (x,y) = 0$ for all $x \neq 0$.
This leads to a singular or badly scaled matrix when using a Newton algorithm to solve
the problem without the presence of a numerical moment.  By adding a numerical moment, the 
resulting system of equations may be better suited for Newton algorithms since 
$\frac{\partial \hF}{\partial P_h^{\pm \mp}} = \frac{\partial F}{\partial P_h^{\pm \mp}} - \alpha$
may be nonsingular even when $P_h^{\pm \mp} \approx 0$.  
For the next numerical test, we let $\alpha = \gamma \mathbf{1}$ for various positive
values of $\gamma$ to see how the numerical moment affects both the accuracy and the performance of
the Newton solver {\it fsolve}.  
The choice for the numerical moment is especially interesting upon noting that 
$\alpha$ is in fact a singular matrix.  
However, with a numerical moment, the perturbation in 
$\frac{\partial \hF}{\partial P_h^{\pm \mp}}$ caused by $P_h^{\pm \mp}$
may be enough to eliminate the singularity since the approximation may now have some curvature.
We let the initial guess be given by the zero function, fix the mesh $N_x = N_y = 20$, and 
let $r=0$.
We can see from Table~\ref{DG_mac0_moment} 
that for $\gamma$ small, {\it fsolve} converges slowly, if at all.  
For $\gamma = 0$, {\it fsolve} does not converge within 100 iterations
even for a very good initial guess.  
However, increasing $\gamma$ does appear to aid {\it fsolve} in its ability to find a root with only 
a small penalty in the approximation error.  
For $r \geq 1$, we again note that Algorithm~\ref{DG_solver_alg} provides a much better suited solver 
due to the degeneracy of the problem.  However, the crux of Algorithm~\ref{DG_solver_alg} reduces 
to a choice of $\gamma > 0$ with $\alpha = \gamma I$ instead of $\alpha = \gamma \mathbf{1}$.
Similar results, as seen in Table~\ref{DG_mac0_moment}, hold for $\alpha = \gamma I$.

\begin{table}[htb]
\caption{Approximation errors when varying $\alpha = \gamma \mathbf{1}$ for Example~\ref{DG_2D_test2} using 
$r=0$, $h$ = 7.071e-02, and  
{\em fsolve} with initial guess $u_h^{(0)} = 0$.
The entry $0^*$ corresponds to an initial guess given by the $L^2$-projection of $u(x,y) = | x |$, 
$q_h^{\pm}(x,y) = \text{sgn } x$, and $P_h^{\mu \nu} (x,y) = 0$ for $\mu, \nu \in \{ + , - \}$.  
The nonlinear solver {\it fsolve} is set to perform a maximum of 100 iterations.}
\begin{center}
\begin{tabular}{| c | c | c | c |}
         \hline
    $\gamma$ & $L^\infty$ norm & $L^2$ norm & {\it fsolve} iterations \\
         \hline
    600 & 2.43e-01 & 2.43e-01 & 9 \\
        \hline
    60 & 2.29e-01 & 2.27e-01 & 9 \\
    	\hline
    12 & 2.02e-01 & 1.98e-01 & 10 \\
    	\hline
    4 & 1.81e-01 & 1.74e-01 & 10 \\
    	\hline
    1 & 3.40e-01 & 2.08e-01 & 100 \\
     	\hline
    $0^*$ & 2.84e-01 & 1.96e-01 & 100 \\
    	\hline
\end{tabular}
\end{center}
\label{DG_mac0_moment}
\end{table}

\begin{example} \label{DG_2D_test3}
Consider the stationary Hamilton-Jacobi-Bellman problem 
\begin{align*}
\min \left\{ - \Delta u , - \Delta u / 2 \right\} &= f  \qquad \textthm{in } \Omega , \\
u &= g \qquad \textthm{on } \partial \Omega, 
\end{align*}
where $\Omega = (0,\pi) \times (-\pi/2,\pi/2)$, 
\[
f(x,y) = \begin{cases}
2 \cos(x) \, \sin(y) , & \textthm{if } (x,y) \in S , \\ 
\cos(x) \, \sin(y) , & \textthm{otherwise} , 
\end{cases}
\]
$S = (0,\pi/2] \times (- \pi/2 , 0] \cup (\pi/2, \pi] \times (0 , \pi/2)$, 
and $g$ is chosen such that
the viscosity solution is given by $u(x,y) = \cos(x) \, \sin(y)$. 
\end{example}

\begin{figure}[htb]
\centering
\includegraphics[scale=0.35]{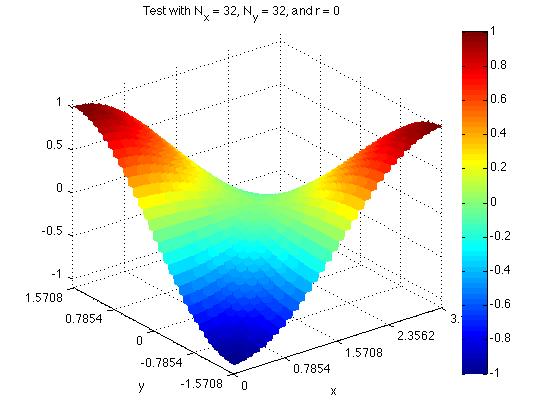}
\caption{
Computed solution for Example~\ref{DG_2D_test3} using 
$r = 0$, $\alpha = 2 I$, $h$ = 1.388e-01, and 
{\em fsolve} with initial guess $u_h^{(0)} = 0$. 
} 
\label{DG_hjb_r0_plot}
\end{figure}

\begin{figure}[htb]
\centering
\includegraphics[scale=0.35]{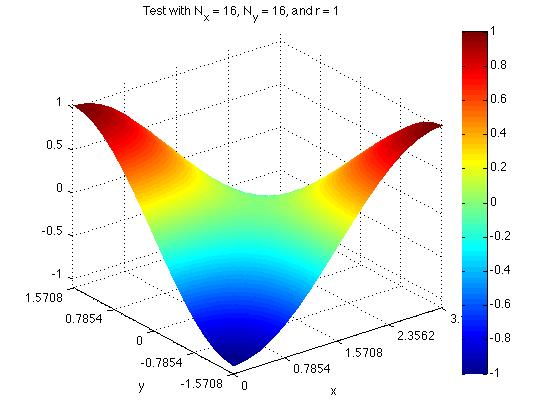}
\caption{
Computed solution for Example~\ref{DG_2D_test3} using 
$r = 1$, $\alpha = 2 I$, $h$ = 2.777e-01, and 
{\em fsolve} with initial guess $u_h^{(0)} = 0$. 
} 
\label{DG_hjb_r1_plot}
\end{figure}

We can see that the optimal coefficient for $\Delta u$ varies over four patches in the domain.
Results for approximating with $r=0,1,2$ can be seen in Tables~\ref{DG_hjb_r0}, 
\ref{DG_hjb_r1}, and \ref{DG_hjb_r2}, respectively,
where we observe optimal convergence rates for $r=0,1$ and near optimal
convergence rates for $r=2$.
Plots for $r=0$ and $r=1$ can be found in Figures~\ref{DG_hjb_r0_plot} and \ref{DG_hjb_r1_plot}.  

\begin{table}[htb]
\caption{
Rates of convergence for Example~\ref{DG_2D_test3} using 
$r = 0$, $\alpha = 2 I$, and {\em fsolve} with initial guess $u_h^{(0)} = 0$. 
}
\begin{center}
\begin{tabular}{| c | c | c | c | c |}
         \hline
    $h$ & $L^\infty$ norm & order & $L^2$ norm & order \\
         \hline
    5.55e-01 & 2.59e-01 &   & 2.73e-01 &   \\
         \hline
    3.70e-01 & 1.63e-01 &  1.14 & 1.75e-01 &  1.10   \\
         \hline
    2.78e-01 & 1.17e-01 &  1.17 & 1.29e-01 &  1.06   \\
         \hline
    1.85e-01 & 7.29e-02 &  1.16 & 8.48e-02 &  1.03   \\
         \hline
    1.39e-01 & 5.27e-02 &  1.13 & 6.33e-02 &  1.02   \\
         \hline
\end{tabular}
\end{center}
\label{DG_hjb_r0}
\end{table}

\begin{table}[htb]
\caption{Rates of convergence for Example~\ref{DG_2D_test3} using 
$r = 1$, $\alpha = 2 I$, and
{\em fsolve} with initial guess $u_h^{(0)} = 0$. 
}
\begin{center}
\begin{tabular}{| c | c | c | c | c |}
         \hline
    $h$ & $L^\infty$ norm & order & $L^2$ norm & order \\
         \hline
    5.55e-01 & 4.89e-02 &   & 2.84e-02 &   \\
         \hline
    3.70e-01 & 2.23e-02 &  1.93 & 1.29e-02 &  1.94   \\
         \hline
    2.78e-01 & 1.27e-02 &  1.97 & 7.38e-03 &  1.95   \\
         \hline
\end{tabular}
\end{center}
\label{DG_hjb_r1}
\end{table}

\begin{table}[htb]
\caption{Rates of convergence for Example~\ref{DG_2D_test3} using 
$r = 2$, $\alpha = 2 I$, and
{\em fsolve} with initial guess $u_h^{(0)} = 0$. 
}
\begin{center}
\begin{tabular}{| c | c | c | c | c |}
         \hline
    $h$ & $L^\infty$ norm & order & $L^2$ norm & order \\
         \hline
    2.22e+00 & 2.82e-01 &   & 1.25e-01 &   \\
         \hline
    7.40e-01 & 9.04e-03 &  3.13 & 9.52e-03 &  2.35   \\
         \hline
    4.44e-01 & 2.39e-03 &  2.60 & 2.88e-03 &  2.34   \\
         \hline
\end{tabular}
\end{center}
\label{DG_hjb_r2}
\end{table}

\begin{example} \label{DG_2D_test4}
Consider the infinite-Laplacian problem 
\begin{align*}
	- \Delta_\infty u := - u_{x x} \, u_{x} \, u_{y}  - u_{x y} \, u_{x} \, u_{y} 
		- u_{y x} \, u_{y} \, u_{y}  - u_{y y} \, u_{y} \, u_{y} &= 0 \qquad \textthm{in } \Omega , \\
	u &=  g \qquad \textthm{on } \partial \Omega,
\end{align*}
where 
$\Omega = (-1,1) \times (-1,1)$ and $g$ is chosen such that
the viscosity solution is given by $u(x,y) =  |x|^{4/3} - |y|^{4/3}$.
While this problem is semilinear and not fully nonlinear, the solution has low 
regularity due to the fact $u \in C^{1,\frac{1}{3}}(\overline{\Omega}) \cap H^1 (\Omega)$.  
\end{example}

\begin{figure}[htb]
\centering
\includegraphics[scale=0.35]{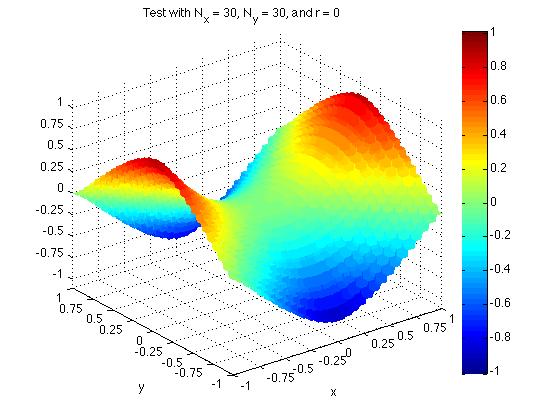}
\caption{
Computed solution for Example~\ref{DG_2D_test4} using 
$r = 0$, $\alpha = 60 I$, $h$ = 9.428e-02, and 
{\em fsolve} with initial guess $u_h^{(0)} = 0$. 
} 
\label{DG_inf_lap_r0_plot}
\end{figure}

\begin{figure}[htb]
\centering
\includegraphics[scale=0.35]{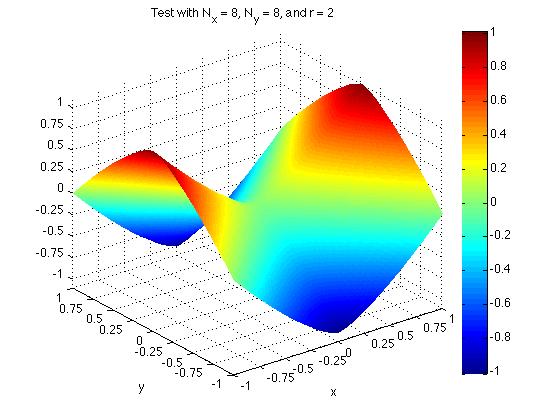}
\caption{
Computed solution for Example~\ref{DG_2D_test4} using 
$r = 2$, $\alpha = 60 I$, $h$ = 3.536e-01, and 
{\em fsolve} with initial guess $u_h^{(0)} = 0$. 
} 
\label{DG_inf_lap_r2_plot}
\end{figure}

By approximation theory, we
expect the error to be bound by $\mathcal{O} (h^{1})$ independent of the degree of the polynomial
basis.  The approximation results for $r=0,1,2$ can be found in 
Tables~\ref{DG_inf_lap_r0}, \ref{DG_inf_lap_r1}, and \ref{DG_inf_lap_r2}, respectively.
Plots for $r=0$ and $r=2$ can be found in 
Figures~\ref{DG_inf_lap_r0_plot} 
and \ref{DG_inf_lap_r2_plot}.
Note that while we observe the theoretical first order bound for the approximation error, 
we also observe that the higher order elements yield more accurate approximations.

\begin{table}[htb]
\caption{Rates of convergence for Example~\ref{DG_2D_test4} using 
$r = 0$, $\alpha = 60 I$, and
{\em fsolve} with initial guess $u_h^{(0)} = 0$. 
}
\begin{center}
\begin{tabular}{| c | c | c | c | c |}
         \hline
    $h$ & $L^\infty$ norm & order & $L^2$ norm & order \\
         \hline
    2.83e-01 & 4.50e-01 &   & 3.37e-01 &   \\
         \hline
    1.41e-01 & 2.83e-01 &  0.67 & 2.02e-01 &  0.74   \\
         \hline
    1.18e-01 & 2.46e-01 &  0.78 & 1.72e-01 &  0.88   \\
         \hline
    9.43e-02 & 2.05e-01 &  0.82 & 1.40e-01 &  0.93   \\
         \hline
\end{tabular}
\end{center}
\label{DG_inf_lap_r0}
\end{table}

\begin{table}[htb] 
\caption{Rates of convergence for Example~\ref{DG_2D_test4} using 
$r = 1$, $\alpha = 60 I$, and
{\em fsolve} with initial guess $u_h^{(0)} = 0$. 
}
\begin{center}
\begin{tabular}{| c | c | c | c | c |}
         \hline
    $h$ & $L^\infty$ norm & order & $L^2$ norm & order \\
         \hline
    4.71e-01 & 4.36e-02 &   & 3.17e-02 &   \\
         \hline
    2.83e-01 & 2.79e-02 &  0.88 & 1.81e-02 &  1.09   \\
         \hline
    2.02e-01 & 2.20e-02 &  0.71 & 1.29e-02 &  1.02   \\
         \hline
\end{tabular}
\end{center}
\label{DG_inf_lap_r1}
\end{table}

\begin{table}[htb] 
\caption{Rates of convergence for Example~\ref{DG_2D_test4} using 
$r = 2$, $\alpha = 60 I$, and
{\em fsolve} with initial guess $u_h^{(0)} = 0$. 
}
\begin{center}
\begin{tabular}{| c | c | c | c | c |}
         \hline
    $h$ & $L^\infty$ norm & order & $L^2$ norm & order \\
         \hline
    5.66e-01 & 2.41e-02 &   & 8.71e-03 &   \\
         \hline
    4.71e-01 & 1.48e-02 &  2.66 & 7.58e-03 &  0.76   \\
         \hline
    3.54e-01 & 1.06e-02 &  1.16 & 4.64e-03 &  1.71   \\
         \hline
\end{tabular}
\end{center}
\label{DG_inf_lap_r2}
\end{table}

\begin{example} \label{DG_time_test5}
Consider the dynamic Hamilton-Jacobi-Bellman problem 
\begin{alignat*}{2}
u_t + \min \left\{ - \Delta u , - \Delta u / 2 \right\} &= f && \qquad \textthm{in } \Omega \times (0,1] , \\
u &= g && \qquad \textthm{on } \partial \Omega \times (0,1], \\  
u &= u_0 && \qquad \textthm{in } \Omega \times \{0\},
\end{alignat*}
where $\Omega = (-1,1) \times (-1,1)$, $f(x,y,t) = s(x,y,t) + 2 \, t \, \left( x \, |x| + y \, |y| \right)$, 
\[
s(x,y,t) = \begin{cases}
 2 t^2 , & \textthm{if } x < 0 \textthm{ and } y < 0 , \\ 
 -4 t^2 , & \textthm{if } x > 0 \textthm{ and } y > 0, \\ 
 0 , & \textthm{otherwise} , 
\end{cases}
\]
and $g$ and $u_0$ are chosen such that
the viscosity solution is given by 
$u(x,y,t) = t^2 \, x \, |x| + t \, y \, |y|$.
Then, for all $t$, we have $u(\cdot, \cdot, t) \in H^2 ( \Omega)$.  
\end{example}

\begin{figure}[htb]
\centering
\includegraphics[scale=0.35]{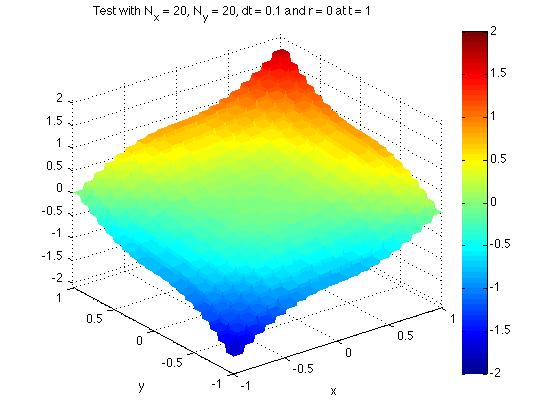}
\caption{
Computed solution at time $t = 1$ for Example~\ref{DG_time_test5} using 
backward Euler time-stepping with 
$r = 0$, $\alpha = 2 I$, $h$ = 1.414e-01, $\Delta t = 0.1$, and 
{\em fsolve} with initial guess $u_h^{0} = \cP_h u_0$. 
} 
\label{DG_hjb_time_r0_plot}
\end{figure}

\begin{figure}[htb]
\centering
\includegraphics[scale=0.35]{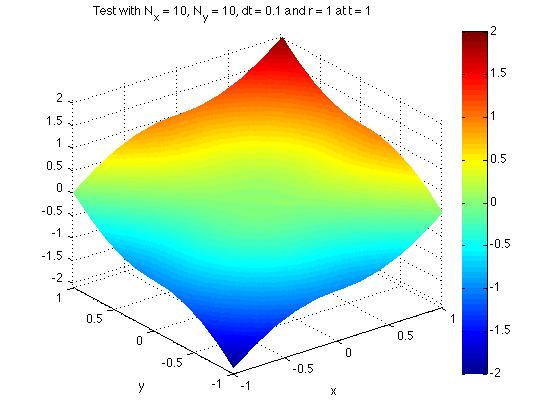}
\caption{
Computed solution at time $t = 1$ for Example~\ref{DG_time_test5} using 
backward Euler time-stepping with 
$r = 1$, $\alpha = 2 I$, $h$ = 2.828e-01, $\Delta t = 0.1$, and 
{\em fsolve} with initial guess $u_h^{0} = \cP_h u_0$. 
} 
\label{DG_hjb_time_r1_plot}
\end{figure}

\begin{figure}[htb]
\centering
\includegraphics[scale=0.35]{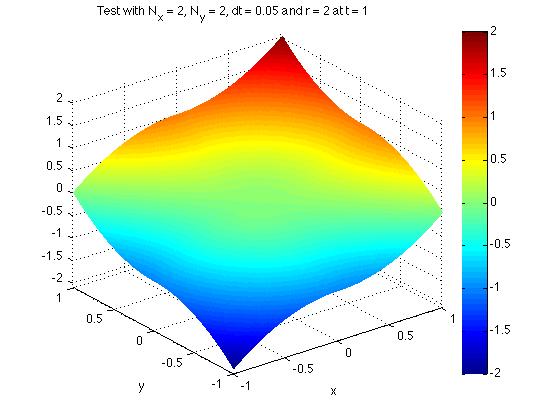}
\caption{
Computed solution at time $t = 1$ for Example~\ref{DG_time_test5} using 
backward Euler time-stepping with 
$r = 2$, $\alpha = 2 I$, $h$ = 1.414, $\Delta t = 0.05$, and 
{\em fsolve} with initial guess $u_h^{0} = \cP_h u_0$. 
} 
\label{DG_hjb_time_r2_plot}
\end{figure}

We expect the spatial rate of convergence to be bound by 2.  However, due to the 
low order time discretization scheme, we can see that our error is dominated by the time 
discretization for $r \geq 1$.  
The spatial orders of convergence for $r=0$ and $r=1$ are recorded in 
Tables~\ref{DG_hjb_time_r0} and \ref{DG_hjb_time_r1}, respectively.
For $r=0$, the spatial discretization order matches the time discretization order, and we do 
observe an optimal rate of convergence.
Using $r=2$, we have the solution $u \in V^h$.  Due to the high level of accuracy when 
using $r=2$, we observe that the time discretization order is in fact 1 as shown in 
Table~\ref{DG_hjb_time_r2}.
Plots for some of the approximations can be found in Figures~\ref{DG_hjb_time_r0_plot}, 
\ref{DG_hjb_time_r1_plot}, and \ref{DG_hjb_time_r2_plot}.

\begin{table}[htb] 
\caption{Rates of convergence in space for Example~\ref{DG_time_test5} 
at time $t=1$ using 
backward Euler time-stepping with 
$r = 0$, $\alpha = 2 I$, $\Delta t = 0.1$, and
{\em fsolve} with initial guess $u_h^{0} = \cP_h u_0$. 
}
\begin{center}
\begin{tabular}{| c | c | c | c | c |}
         \hline
    $h$ & $L^\infty$ norm & order & $L^2$ norm & order \\
         \hline
    2.83e-01 & 5.62e-01 &   & 2.63e-01 &   \\
         \hline
    1.77e-01 & 3.62e-01 &  0.93 & 1.71e-01 &  0.92   \\
         \hline
    1.41e-01 & 2.92e-01 &  0.96 & 1.38e-01 &  0.96   \\
         \hline
\end{tabular}
\end{center}
\label{DG_hjb_time_r0}
\end{table}

\begin{table}[htb] 
\caption{Rates of convergence in space for Example~\ref{DG_time_test5} 
at time $t=1$ using 
backward Euler time-stepping with 
$r = 1$, $\alpha = 2 I$, $\Delta t = 0.1$, and
{\em fsolve} with initial guess $u_h^{0} = \cP_h u_0$. 
}
\begin{center}
\begin{tabular}{| c | c | c | c | c |}
         \hline
    $h$ & $L^\infty$ norm & order & $L^2$ norm & order \\
         \hline
    4.71e-01 & 7.41e-02 &   & 5.00e-02 &   \\
         \hline
    3.54e-01 & 4.21e-02 &  1.96 & 3.56e-02 &  1.18   \\
         \hline
    2.83e-01 & 3.10e-02 &  1.38 & 2.76e-02 &  1.14   \\
         \hline
\end{tabular}
\end{center}
\label{DG_hjb_time_r1}
\end{table}

\begin{table}[htb] 
\caption{Rates of convergence in time for Example~\ref{DG_time_test5} 
at time $t=1$ using 
backward Euler time-stepping with 
$r = 2$, $\alpha = 2 I$, $h$ = 1.414, and
{\em fsolve} with initial guess $u_h^{0} = \cP_h u_0$. 
}
\begin{center}
\begin{tabular}{| c | c | c | c | c |}
         \hline
    $\Delta t$ & $L^\infty$ norm & order & $L^2$ norm & order \\
         \hline
    5.00e-01 & 4.12e-02 &   & 4.12e-02 &   \\
         \hline
    2.50e-01 & 2.11e-02 &  0.96 & 2.11e-02 &  0.97   \\
         \hline
    1.00e-01 & 8.55e-03 &  0.99 & 8.49e-03 &  0.99   \\
         \hline
    5.00e-02 & 4.29e-03 &  1.00 & 4.25e-03 &  1.00   \\
         \hline
\end{tabular}
\end{center}
\label{DG_hjb_time_r2}
\end{table}

\section{Conclusion} \label{sec-6}

In this paper, we have formulated a framework for designing LDG methods that 
approximate the viscosity solution of fully nonlinear second order elliptic and parabolic 
PDEs in high dimensions.  
We then focused on a particular LDG method within the framework 
that corresponded to the Lax-Friedrichs-like numerical operator.   
The key tools in designing the numerical operator are the introduction of 
a numerical viscosity and numerical moment.  
Through numerical tests, we observed the potential for the given framework 
that was originally motivated by successful numerical techniques for Hamilton-Jacobi equations 
as well a FD framework that abstracts the indirect techniques of the vanishing moment method.  

A major task when approximating viscosity solutions is designing methods 
that are selective enough to rule out low regularity artifacts based on the PDE 
yet flexible enough to account for the fact that the viscosity solution itself may 
have low regularity.  
Numerical tests in this paper as well as \cite{Feng_Lewis12a,Feng_Lewis12c,Feng_Lewis13} 
indicate that the dichotomy between successfully capturing a smooth solution 
by ruling out lower regularity artifacts while still being able to approximate low regularity 
functions is an issue that may best be tackled when discretizing a PDE and designing a solver 
for the resulting algebraic system occur in concert with each other.  
Our numerical tests show that the numerical moment successfully removes numerical artifacts 
in many examples.  
However, they also indicate that the numerical moment alone cannot rule out 
all numerical artifacts in all instances.  
The best hope is that the discretization can effectively destabilize numerical artifacts when 
paired with an appropriate solver, 
as was achieved in all of our test problems that used the inverse-Poisson fixed-point solver.  
Given the observed potential for destabilizing low-regularity PDE artifacts using a numerical moment, 
another promising direction of research is using the numerical moment as a low-regularity 
indicator when designing and implementing adaptive methods.

\bibliographystyle{plain}

\end{document}